\documentclass[12pt]{amsart}
\usepackage{amssymb}
\usepackage{amsfonts}
\usepackage{amsmath}
\usepackage{array}
\usepackage[ps,matrix,curve]{xypic}
\usepackage{pstricks}

\def\SetTableau#1#2#3#4{%
  \gdef\Tabvrule{\vrule\vrule width-0.4pt}
  \gdef\Tabhrule{\hrule\hrule height-0.4pt}  
  \gdef\Tabstrut{\vrule height#1 depth#2 width0pt\relax}
  \gdef\Tabbox##1{\hbox to #3{\hskip0.4pt\hfill\Tabstrut$#4##1$\hfill}}
} 
\def\TasseTableau{\SetTableau{1.48ex}{0.32ex}{1.8ex}{\scriptstyle}}
\def\PetitTableau{\SetTableau{1.65ex}{0.55ex}{2.2ex}{\scriptstyle}}

\def\Case#1{\vcenter{\Tabhrule%
                   \hbox{\Tabvrule\Tabbox{#1}\Tabvrule}\Tabhrule}}

\def\GenTab#1{\vcenter{\halign{&$\Case{##}$\cr#1}}\egroup}
\def\Tableau#1{%
  \bgroup%
  \let\ =\omit%
  \let\\=\cr%
  \offinterlineskip\GenTab{#1}}
\let\ifdraft=\iffalse
\let\ifdraft=\iftrue

\def\QMR{{\it QMR}}
\def\MR{{\bf MR}}

\def\Chi{\operatorname{\rm X}}
\def\concat{{\,}}

\def\AKS{\operatorname{\mathcal H}}
\def\AKSq{\AKS_{n,r}(q)}
\def\AKSz{\AKS_{n,r}(0)}
\def\AKSzgen#1#2{\AKS_{#1,#2}(0)}

\newtheorem{example}{Example}[section]
\newtheorem{note}[example]{Note}
\newtheorem{theorem}[example]{Theorem}

\newtheorem{corollary}[example]{Corollary}

\newtheorem{proposition}[example]{Proposition}
\newtheorem{algorithm}[example]{Algorithm}
\newtheorem{lemma}[example]{Lemma}

\def\qed{\hspace{3.5mm} \hfill \vbox{\hrule height 3pt depth 2 pt width 2mm}
\bigskip}

\def\FQSym{{\bf FQSym}}

\def\NCSF{{\bf Sym}}
\def\QSym{{\it QSym}}

\def\ssh{\Cup}

\def\sconc{\bullet}
\def\Std{{\rm Std}}
\def\rsw{{\rm RStd}}

\def\convol{{*}}

\def\<{\langle}
\def\>{\rangle}

\def\CC{\operatorname{\mathbb C}}

\def\ZZ{\operatorname{\mathbb Z}}
\def\NN{\operatorname{\mathbb N}}

\def\F{{\bf F}}

\def\G{{\bf G}}

\def\SG{{\mathfrak S}}

\def\goth{\mathfrak}

\def\Des{\operatorname{C}}
\def\Dess{\operatorname{D}}

\def\dim{{\rm dim}}

\def\shuff#1#2{\mathbin{
\hbox{\vbox{ \hbox{\vrule \hskip#2 \vrule height#1 width 0pt
}%
\hrule}%
\vbox{ \hbox{\vrule \hskip#2 \vrule height#1 width 0pt
\vrule }%
\hrule}%
}}}

\def\shuf{{\mathchoice{\shuff{7pt}{3.5pt}}%
{\shuff{6pt}{3pt}}%
{\shuff{4pt}{2pt}}%
{\shuff{3pt}{1.5pt}}}}%
\def\shuffle{\,\shuf\,}

\def\finer{\leq}

\def\x{{\bf x}}
\def\c{{\bf c}}
\def\d{{\bf d}}
\def\rlambda{{\boldsymbol \lambda}}

\def\antiinvol{\varphi}
\def\bijcyc{\phi}
\def\s{s} 

\def\rad{\operatorname{rad}}

\def\Ext{\operatorname{Ext}}

\title[Yang-Baxter bases of $H_n(0)$ and representation theory of $\AKSz$]
{Yang-Baxter bases of $0$-Hecke algebras and representation theory of\\
$0$-Ariki-Koike-Shoji algebras}

\author[F. Hivert, J.-C.~Novelli, and J.-Y.~Thibon]%
{Florent Hivert, Jean-Christophe Novelli, and Jean-Yves Thibon}

\address[] {Institut Gaspard Monge, Universit\'e de Marne-la-Vall\'ee \\
5 Boulevard Descartes \\Champs-sur-Marne \\77454 Marne-la-Vall\'ee cedex 2 \\
FRANCE}
\email[Florent Hivert]{hivert@univ-mlv.fr}
\email[Jean-Christophe Novelli]{novelli@univ-mlv.fr}
\email[Jean-Yves Thibon]{jyt@univ-mlv.fr} 
\date{}

\begin{document}

\begin{abstract}
After reformulating the representation theory of $0$-Hecke algebras in an
appropriate family of Yang-Baxter bases, we investigate certain
specializations of the Ariki-Koike algebras, obtained by setting $q=0$ in a
suitably normalized version of Shoji's presentation.
We classify the simple and projective modules, and describe restrictions,
induction products, Cartan invariants and decomposition matrices.
This allows us to identify the Grothendieck rings of the towers of algebras in
terms of certain graded Hopf algebras known as the Mantaci-Reutenauer descent
algebras, and Poirier quasi-symmetric functions.
We also describe the Ext-quivers, and conclude with numerical tables.
\end{abstract}

\maketitle

\tableofcontents
\section{Introduction}

Given an \emph{inductive tower of algebras}, that is, a sequence of
algebras
\begin{equation}
A_0 \hookrightarrow A_1 \hookrightarrow \ldots \hookrightarrow A_n
\hookrightarrow \ldots,
\end{equation}
with embeddings $A_m\otimes A_n \hookrightarrow A_{m+n}$ satisfying
appropriate compatibility conditions, one can introduce two \emph{Grothendieck
rings}
\begin{equation}
{\mathcal G}:=\bigoplus_{n\ge 0}G_0(A_n)\,,\quad
{\mathcal K}:=\bigoplus_{n\ge 0}K_0(A_n)\,,
\end{equation}
where $G_0(A)$ and $K_0(A)$ are the (complexified) Grothendieck groups of the
categories of finite-dimensional $A$-modules and projective $A$-modules
respectively, with multiplication of the classes of an $A_m$-module $M$ and an
$A_n$-module $N$ defined by
\begin{equation}
[M] \cdot [N] = [M\widehat{\otimes} N] =
[M\otimes N \uparrow_{A_m\otimes A_n}^{A_{m+n}}].
\end{equation}

On each of these Grothendieck rings, one can define a coproduct by means of
restriction of representations. Under favorable circumstances, this turns
these two rings into mutually dual Hopf algebras.

The basic example of this situation is the character ring of symmetric groups
(over~$\CC$), due to Frobenius. Here the $A_n=\CC\SG_n$ are semi-simple
algebras, so that
\begin{equation}
{G}_0(A_n) = {K}_0(A_n)= R(A_n),
\end{equation}
where $R(A)$ denotes the vector space spanned by isomorphism classes of
indecomposable modules which in this case are all simple and projective.
The irreducible representations $[\lambda]$ of $A_n$ are parametrized by
partitions $\lambda$ of $n$, and the Grothendieck ring is isomorphic to the
algebra $Sym$ of symmetric functions under the
correspondence
$[\lambda] \leftrightarrow s_\lambda$,
where $s_\lambda$ denotes the Schur function associated with $\lambda$.

Other known examples with towers of group algebras over the complex numbers,
$A_n=\CC G_n$, include the cases of wreath products
$G_n = \Gamma\wr\SG_n$ (Specht), finite linear groups
$G_n = GL(n,\F_q)$ (Green), \emph{etc.}, all related to symmetric functions
(see~\cite{Mcd,Zel}).

Examples involving modular representations of finite groups and non-semisimple
specializations of Hecke algebras have also been worked out
(see~\cite{BK,LLT,LT}).
For example, finite Hecke algebras of type $A$ at roots of unity
($A_n=H_n(\zeta)$, $\zeta^k=1$) yield quotients and subalgebras of $Sym$
\cite{LLT}
\begin{equation}
{\mathcal G} = Sym/(p_{km}=0), \quad
{\mathcal K} = \CC\left[\,p_i\,|\,i\not\equiv0\pmod k\right],
\end{equation}
where the $p_i$ are the power sum symmetric functions, supporting level $1$
irreducible representations of the affine Lie algebra $\widehat{\goth{sl}}_k$,
while Ariki-Koike algebras at roots of unity give rise to higher level
representations of the same Lie algebras \cite{Ari}.
The $0$-Hecke algebras $A_n=H_n(0)$ correspond to the pair Quasi-symmetric
functions - Noncommutative symmetric functions, ${\mathcal G}=\QSym$,
${\mathcal K}=\NCSF$ \cite{NCSF4}.
Affine Hecke algebras at roots of unity lead to $U^+(\widehat{\goth{sl}}_k)$
and $U^+(\widehat{\goth{sl}}_k)^*$ \cite{Ari}, and the cases of generic affine
Hecke algebras can be reduced to a subcategory admitting as Grothendieck
rings $U^+(\widehat{\goth{gl}}_\infty)$ and
$U^+(\widehat{\goth{gl}}_\infty)^*$ \cite{Ari}.

A further interesting example is the tower of $0$-Hecke-Clifford algebras
\cite{Ols,BHT}, giving rise to the peak algebras \cite{NCSF2,Stem}.

Here, we shall show that appropriate versions at $q=0$ of the Ariki-Koike
algebras derived from the presentation of Shoji \cite{Sho,Sasho} admit as
Grothendieck rings two known combinatorial Hopf algebras, the
Mantaci-Reutenauer descent algebras (associated with the corresponding wreath
products) \cite{MR}, and their duals, a generalization of quasi-symmetric
functions, introduced by Poirier in~\cite{Poi} and more recently considered
in~\cite{NT-coul,BH}.
This work can be understood as the solution of an \emph{inverse problem} in
representation theory: given the Grothendieck rings, reconstruct the algebras.
The main point is that the standard presentations of Ariki-Koike algebras at
$q=0$ do not give the required Grothendieck rings.

\bigskip
This article is structured as follows. After recalling some classical
definitions, we introduce new combinatorial structures (cycloribbons and
anticycloribbons) needed in the sequel (Section 2). Next, we study a limit
case of the Yang-Baxter bases of Hecke algebras introduced in~\cite{LLT:Flag},
which will be one of our main tools (Section 3).
In Section 4, we finally introduce the fundamental objects of our study,
the $0$-Ariki-Koike-Shoji algebras $\AKSz$, and special bases, well suited
for analyzing representations.
Next, we obtain the classification of simple $\AKSz$-modules, which turn out
to be all one-dimensional, labelled by the $r(r+1)^{n-1}$ cycloribbons.
We then describe induction products and restrictions of these simple modules,
which allows us to identify the first Grothendieck ring $\mathcal G$ with a
Hopf subalgebra of Poirier's Quasi-symmetric functions, dual to the
Mantaci-Reutenauer Hopf algebra (Section 5).
Since $\AKSz$ is self-injective, this duality gives immediately the
Grothendieck ring $\mathcal K$ associated with projective modules. An
alternative labelling of the indecomposable projective modules leads then to a
simple description of basic operations such as induction products, restriction
to $H_n(0)$, or induction from a projective $H_n(0)$-module. Summarizing, we
obtain an explicit description of the Cartan-Brauer triangle, in particular of
the Cartan invariants and of the decomposition matrices (Section 6).
We conclude with a description of the Ext-quiver of $\AKSz$. The last section
contains tables of $q$-Cartan invariants (recording the radical filtrations of
projective modules) and decomposition matrices.

\section{Notations and preliminaries}

\subsection{Words and permutations}

Let $w$ be a word on a totally ordered alphabet $A$.
We denote by $\overline{w}$ the \emph{mirror image} (reading from right to
left) of $w$.
If the alphabet is $\NN^*$, the \emph{evaluation} of $w$ is the sequence of
number of occurrences of $1$, $2$, and so on until all letters of $w$ have
been seen.
For example, the evaluation of $15423341511457$ is $(4,1,2,3,3,0,1)$.

We will generally use the operation of concatenation on words but we will also
use the \emph{shifted concatenation} $u\sconc v$ of two words $u$ and $v$ over
the positive integers consisting in concatenating $u$ with
$v[k]:= (v_1+k,\ldots,v_l+k)$, where $k=|u|$ is the length of $u$.
For example $1431\sconc232 = 1431676$.

We shall also make use of another classical operation on words, known as
\emph{standardization}. The standardized word $\Std(w)$ of $w$ is the
permutation having the same inversions as $w$,
and its right standardized word is
$\rsw(w) := \overline{\Std(\overline{w})}$.
For example, $\Std(412461415) = 514692738$,
$\rsw(412461415) = 734692518$.

We will represent permutations as words of size $n$ over the alphabet
$\{1,\ldots,n\}$.
Denote by $s_i$ the transposition exchanging $i$ and $i+1$, and by
$\omega_n$ the maximal permutation $(n,n-1,\ldots,1)$. For
$\sigma=(\sigma_1,\ldots,\sigma_n)$, we set
\begin{equation}
\overline{\sigma} := (\sigma_n,\ldots,\sigma_1)=\sigma \omega_n,\quad
\sigma^\# := (n+1-\sigma_1,\ldots,n+1-\sigma_n) = \omega_n \sigma.
\end{equation}

\subsection{Compositions}

Recall that a \emph{composition} of an integer $n$ is any finite sequence of
positive integers $I=(i_1,\ldots,i_k)$ of sum $|I|:=i_1+\cdots+i_k=n$. It can
be pictured as a ribbon diagram, that is, a set of rows composed of square
cells of respective lengths $i_j$, the first cell of each row being attached
under the last cell of the previous one. $I$ is called the \emph{shape} of the
ribbon diagram. The \emph{conjugate} $I^{\sim}$ of a composition $I$ is
the composition built with the number of cells of the columns of its ribbon
diagram read from \emph{right to left}. Its \emph{mirror image} $\overline{I}$
is the reading of $I$ from right to left.
For example, if $I=(2,3,1,2)$, $I^\sim$ is $(1,3,1,2,1)$ and
$\overline{I}$ is $(2,1,3,2)$.
\begin{figure}[ht]
\begin{equation*}
\PetitTableau
I =
\Tableau{ & \\\ & & & \\\ &\ &\ & \\\ &\ &\ & & \\}
\quad
I^{\sim} =
\Tableau{\\ & & \\\ &\ &\\\ &\ & &\\\ &\ &\ & \\}
\quad
\overline{I} = 
\Tableau{ & \\\ &\\\ & & & \\\ &\ &\ && \\}
\end{equation*}
\caption{Ribbon diagrams of the compositions
$I=(2,3,1,2)$, $I^\sim=(1,3,1,2,1)$, and $\overline{I}=(2,1,3,2)$.}
\end{figure}

Given a filling of a ribbon diagram, we define its \emph{row reading} as the
reading of its rows from left to right and from top to bottom.
Finally, the \emph{descent set} $\Dess(I)$ of a composition
$I=(i_1,\ldots,i_l)$ is the set $\{i_1,i_1+i_2,\ldots,i_1+\cdots+i_{l-1}\}$
composed of the \emph{descents} of $I$.
Recall also that the \emph{descent set} $\Dess(\sigma)$ of a permutation
$\sigma$ is the set of $i$ such that $\sigma(i)>\sigma(i+1)$ (the
\emph{descents} of $\sigma$), and the \emph{descent composition}
$\Des(\sigma)$ of $\sigma$ is the unique composition $I$ of $n$ such that
$\Dess(I)=\Dess(\sigma)$, that is, the shape of a filled ribbon diagram whose
row reading is $\sigma$ and whose rows are increasing and columns decreasing.
For example, Figure~\ref{perm-compo} shows that the descent composition of
$(3,5,4,1,2,7,6)$ is $I=(2,1,3,1)$.

\begin{figure}[ht]
$
\PetitTableau
\Tableau{3 &5 \\\ & 4\\\ & 1& 2& 7 \\\ &\ &\ & 6 \\ }
$
\caption{\label{perm-compo}The ribbon diagram of the permutation
$3541276$.}
\end{figure}

Conversely, with a composition $I$, associate its \emph{maximal permutation}
$\sigma=\omega(I)$ as the permutation with descent composition
$I$ and maximal inversion number.
Similarly, the minimal permutation $\alpha(I)$ is the permutation with descent
composition $I$ and minimal inversion number.
For example, if $I=(2,1,3,1)$, $\omega(I)=6752341$ and
$\alpha(I)=1432576$.

\subsection{Shuffle, shifted shuffle, and convolution}

The \emph{shuffle product} $\shuffle$ of two words $u$ and $v$ over an
alphabet $A$ is inductively defined by:
\begin{equation}
\label{shuffle}
  \text{if $u=au'$ and $v=bv'$ with $a,b\in A$,\quad then \quad}
  u\shuffle v = a(u' \shuffle v) +  b(u\shuffle v'),
\end{equation}
with the initial conditions $u\shuffle \varepsilon=\varepsilon\shuffle u =u $,
$\varepsilon$ being the empty word.

The \emph{shifted shuffle} $\ssh$ of two permutations
$\sigma'\in\SG_k$ and $\sigma''\in\SG_l$ is the shuffle 
of $\sigma'$ and $\sigma''[k]$. 
For example,
\begin{equation}
21\ssh 12    = 2134 + 2314 + 3214 + 2341 + 3241 + 3421.
\end{equation}

The \emph{convolution product} $\convol$ of two permutations
$\sigma'\in\SG_k$ and $\sigma''\in\SG_l$ is the sum of permutations $\sigma$
in $\SG_{k+l}$ such that $\Std(\sigma_1\cdots\sigma_k)=\sigma'$ and 
$\Std(\sigma_{k+1}\cdots\sigma_{k+l})=\sigma''$. Equivalently, it is the sum
obtained by inverting the permutations occuring in the shifted shuffle of the
inverses of $\sigma'$ and $\sigma''$.
It is well-known that the convolution of two permutations is an interval
of the (left) weak order. For example, one has
\begin{equation}
21\convol 12 = 2134 + 3124 + 3214 + 4123 + 4213 + 4312 = [2134,4312].
\end{equation}

In all the paper, we make use of a set $C=\{1,\ldots,r\}$, called the
\emph{color set}. The {color words} will be written in bold type.

The previous definitions can be extended to colored permutations as
in~\cite{NT-coul}, that is, pairs $(\sigma,\c)\in \SG_n\times C^n$.
The \emph{shifted shuffle} of two such elements
$(\sigma',\c')\ssh(\sigma'',\c'')$ is the sum of colored permutations obtained
by shuffling with shift the two permutations, the colors remaining attached to
letters.
For example, representing colors as exponents,
\begin{equation}
\begin{split}
2^21^1\ssh 1^12^2 &= 2^21^1\shuffle 3^14^2 \\
&= 2^21^13^14^2 + 2^23^11^14^2 + 3^12^21^14^2
+ 2^23^14^21^1 + 3^12^24^21^1 + 3^14^22^21^1.
\end{split}
\end{equation}

The \emph{convolution} $(\sigma',\c')\convol(\sigma'',\c'')$ of two colored
permutations is the sum of colored permutations $(\sigma,\c'\concat\c'')$,
where $\sigma$ runs over $\sigma'\convol\sigma''$. Equivalently, if one
defines the \emph{inverse} of a colored permutation $(\sigma,\c)$ as
$(\sigma^{-1},\c\cdot\sigma^{-1})$,
the convolution of two colored permutations is the sum of the inverses of
the elements occuring in the shifted shuffle of their inverses.
For example,
\begin{equation}
2^21^1\convol 1^12^2 = 2^21^13^14^2 + 3^21^12^14^2 + 3^22^11^14^2
+ 4^21^12^13^2 + 4^22^11^13^2 + 4^23^11^12^2.
\end{equation}

The sums being multiplicity free, we can regard these as sets, and write, 
\emph{e.g.}, $2314 \in 21\ssh 12$, or $3^21^12^14^2\in2^21^1\convol 1^12^2$.

\subsection{Cycloribbons and anticycloribbons}

Let $I$ be a composition of $n$ and $\c\in C^n$ be a color word of
length~$n$. The pair $[I,\c]$ is called a \emph{colored ribbon},
and depicted as the filling of $I$ whose row reading is $\c$.
We say that this filling is a \emph{cyclotomic ribbon} (\emph{cycloribbon}
for short) if it is \emph{weakly increasing} in rows and
\emph{weakly decreasing} in columns.
Notice that there are $r(r+1)^{n-1}$ cycloribbons with at most $r$ colors
since when building the ribbon cell by cell, one has $r$ possibilities for its
first cell and then $r+1$ possibilities for the next ones: $1$ possibility for
the $r-1$ choices different from the previous one and $2$ possibilities for
the same choice (right or down). Here are the five cycloribbons of shape
$(2,1)$ with two colors:
\begin{equation}
\label{5cyclo}
\PetitTableau
\Tableau{1&1\\\ &1\\ }
\qquad
\Tableau{1&2\\\ &1\\ }
\qquad
\Tableau{1&2\\\ &2\\ }
\qquad
\Tableau{2&2\\\ &1\\ }
\qquad
\Tableau{2&2\\\ &2\\ }
\end{equation}

We say that a colored ribbon is an \emph{anticyclotomic ribbon}
(\emph{anticycloribbon} for short) if it is \emph{weakly decreasing} in rows
and \emph{weakly increasing} in columns. There are as many anticycloribbons
as cycloribbons. The relevant bijection $\bijcyc$ from one set to the other is
the restriction of an involution on all colored ribbons: read a ribbon $R$
row-wise and build its image $\bijcyc(R)$ cell by cell as follows:
\begin{itemize}
\item if the $(i+1)$-th cell has the same content as the $i$-th cell, glue it
at the same position as in $R$ (right or down),
\item if the $(i+1)$-th cell does not have the same content as the $i$-th
cell, glue it at the other position (right or down).
\end{itemize}

\noindent
For example, the colored ribbons of Figure~\ref{cycanti} are exchanged by
$\bijcyc$:

\begin{figure}[ht]
\begin{equation*}
\PetitTableau
\Tableau{1&1&3\\\ &\ & 3\\\ &\ &3\\\ &\ &2\\\ &\ &1&1&4&5\\ }
\qquad \qquad \qquad
\Tableau{1&1\\\ &3\\\ &3\\\ &3&2&1&1\\\ &\ &\ &\ &4\\\ &\ &\ &\ &5\\ }
\end{equation*}
\caption{\label{cycanti}A cycloribbon $R$ and the anticycloribbon
$\bijcyc(R)$.}
\end{figure}

Let us now extend to colored permutations and anticycloribbons the
correspondances between permutations and ribbons:
we say that $i$ is an \emph{anti-descent} of $(\sigma,\c)$ iff
$c_i<c_{i+1}$ or, $c_i=c_{i+1}$ and $\sigma_i<\sigma_{i+1}$.
Given a colored permutation $(\sigma,\c)$, define its associated
anticycloribbon as the pair $[I,\c]$ where $I$ has a descent at the
$i$-th position iff $i$ is an anti-descent of the colored permutation.
Conversely, with a given anticycloribbon $[I,\c]$, associate its
\emph{maximal colored permutation} ${(\sigma,\c)}$ where $\sigma=\omega(J)$,
where $\Dess(J)= \{
i\,|\, \text{$c_i\not= c_{i+1}$, or $c_i=c_{i+1}$ and $i\not\in\Dess(I)$}\}$.


For example, the anti-descents of $(81732564,11132221)$ are $\{2,3,5,6\}$.
Its associated anticycloribbon is
\begin{equation}
\PetitTableau
\Tableau{1&1\\\ & 1\\\ & 3 &2\\\ &\ &2 \\\ &\ &2 &1\\ }
\end{equation}
and the maximal colored permutation having this associated anti-cycloribbon
is $\omega(1,2,1,3,1)= (86752341,11132221)$.

\newpage
\section{The Hecke algebra of type $A$ at $q=0$ and Yang-Baxter bases}

In this section, we first recall some combinatorial aspects of the
representation theory of $H_n(0)$, such as a realization of its simple and
indecomposable projective modules in the left regular representation and the
calculation of the composition factors of an induction product of two simple
or indecomposable projective modules.
We will show in the next Section how these constructions can be generalized to
the $0$-Ariki-Koike-Shoji algebra $\AKSz$.
To achieve this, we will need to investigate in some detail a special case of
the so-called Yang-Baxter bases of $H_n(0)$.

\subsection{Representation theory}
\label{rthn}

Let us consider the Hecke algebra $H_n(0)$ of type $A$ at $q=0$ in the
presentation of~\cite{Nor,NCSF4}, that is, the $\CC$-algebra with
generators $T_i$, $1\leq i\leq n-1$ and relations:
\begin{alignat}{2}
T_i(1+T_i)      & = 0                   && (1\leq i\leq n-1), \\
T_i T_{i+1} T_i & = T_{i+1} T_i T_{i+1} &\quad& (1\leq i\leq n-2), \\
T_i T_j         & = T_j T_i             && (|i-j|\geq2).
\end{alignat}

Let $\sigma =: \sigma_{i_1}\dots\sigma_{i_p}$ be a reduced word
for a permutation $\sigma\in\SG_n$. The defining relations of $H_n(0)$ ensure
that the element $T_{\sigma} := T_{i_1}\dots T_{i_p}$ is independent of the
chosen reduced word for $\sigma$. Moreover, the well-defined family
$(T_{\sigma})_{\sigma\in\SG_n}$ is a basis of the Hecke algebra, which is
consequently of dimension $n!$.

\subsubsection{Simple modules and induction}

It is known that $H_n(0)$ has $2^{n-1}$ simple modules, all one-dimensional,
naturally labelled by compositions $I$ of $n$~\cite{Nor}:
following the notation of~\cite{NCSF4}, let $\eta_I$ be the generator of the
simple $H_n(0)$-module $S_I$ associated with $I$ in the left regular
representation. It satisfies
\begin{equation}
\label{etaI}
\left\{
\begin{array}{r@{\,}cr}
-T_i \eta_I   &= \eta_I & \text{if $i\in\Dess(I)$,}\\
(1+T_i)\eta_I &= \eta_I & \text{otherwise.}
\end{array}
\right.
\end{equation}

The composition factors of the induction product $S_I\widehat\otimes S_J$
of two simple $0$-Hecke modules are easily described in terms of
permutations.

\begin{figure}[ht]
\PetitTableau
\ifdraft
\rotateleft{
\newdimen\vcadre\vcadre=0.2cm 
\newdimen\hcadre\hcadre=0.0cm 
\def\GrTeXBox#1{\vbox{\vskip\vcadre\hbox{\hskip\hcadre%
      $\Tableau{#1\\ }$%
   \hskip\hcadre}\vskip\vcadre}}
\def\GrTeXBoxb#1{\vbox{\vskip\vcadre\hbox{\hskip\hcadre%
   $#1$\hskip\hcadre}\vskip\vcadre}}
  \def\loopRight{\ar @`{[]+<1.2cm,1.1cm>,[]+<1.2cm,-1.1cm>}}
  \def\loopLeft{\ar @`{[]+<-1.2cm,1.1cm>,[]+<-1.2cm,-1.1cm>}}

$\vcenter{\xymatrix@R=2.2cm@C=+0.5cm{
   & *{\GrTeXBoxb{1 2 3 4}} \ar @{->}[d]^{2} & \\
   & *{\GrTeXBoxb{1 3 2 4}}
   \ar @{..>}[dl]_{1} \ar @{=>}[dr]^{3}& \\
 *{\GrTeXBoxb{2 3 1 4}}\ar @{=>}[dr]_{3}
 & & *{\GrTeXBoxb{1 4 2 3}} \ar @{..>}[dl]^{1}\\
   & *{\GrTeXBoxb{2 4 1 3}} \ar @{->}[d]^{2}& \\
   & *{\GrTeXBoxb{3 4 1 2}} & \\
}}
\quad\longmapsto\quad
\vcenter{\xymatrix@R=2.2cm@C=+0.4cm{
   & *{\GrTeXBoxb{b_{1234}}} \loopRight^{{}^{1+T_1}_{1+T_3}}\ar @{->}[d]^{T_2} & \\
   & *{\GrTeXBoxb{b_{1324}}}
   \ar @{..>}[dl]_{T_1} \loopRight^{{}^{-T_2}}\ar @{=>}[dr]^{T_3}& \\
 *{\GrTeXBoxb{b_{2314}}}\loopLeft_{{}^{-T_1}_{1+T_2}}\ar @{=>}[dr]_{T_3}
 & & *{\GrTeXBoxb{b_{1423}}}\loopRight^{{}^{1+T_2}_{-T_3}} \ar @{..>}[dl]^{T_1}\\
   & *{\GrTeXBoxb{b_{2413}}}\loopRight^{{}^{-T_1}_{-T_3}} \ar @{->}[d]^{T_2}& \\
   & *{\GrTeXBoxb{b_{3412}}}\loopRight^{{}^{1+T_1}_{-T_2,1+T_3}} & \\
}}
\quad\longmapsto\quad
\vcenter{\xymatrix@R=1.6cm@C=+0.2cm{
   & *{\GrTeXBox{ &  &  & }} \ar @{->}[d] & \\
   & *{\GrTeXBox{& \\\ & & }}
   \ar @{->}[dl] \ar @{->}[dr]& \\
 *{\GrTeXBox{\\  &  & }}\ar @{->}[dr]
 & & *{\GrTeXBox{&  & \\\ &\ & }}\ar @{->}[dl]\\
   & *{\GrTeXBox{\\ & \\\ &}}\ar @{->}[d]& \\
   & *{\GrTeXBox{& \\\ & &}} & \\
}}
$
}
\fi
\caption{\label{0-Hecke-graphe}%
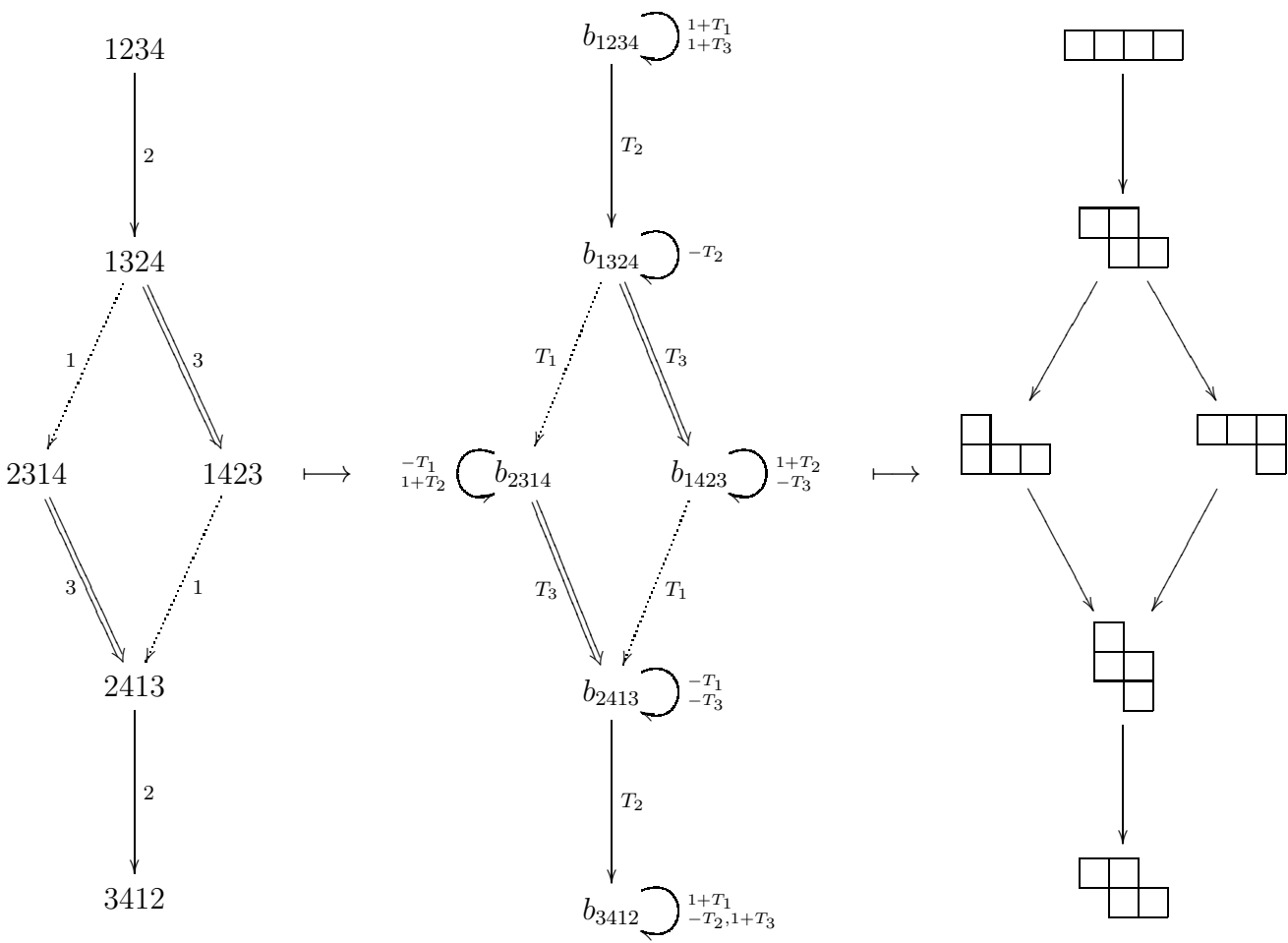
Induction product of two simple modules of $H_2(0)$ to $H_4(0)$. The first
graph is the restriction of the left weak order to the convolution of
$12\convol12$.
The second graph is its $0$-Hecke graph and the third graph represents the
composition factors, read from the descent compositions of the inverse
permutations.}
\end{figure}

\medskip
Let us say that a family $(b)$ of elements of a $H_n(0)$-module $M$ is a
\emph{combinatorial basis} if it is a basis of $M$ such that $T_i\,b$ is
either $-b$, or $0$ or another basis element $b'$ for all $i$. Then $M$ is
completely encoded by the edge-labelled directed graph having as vertices the
basis elements, with a loop labelled $-T_i$ iff $-T_i\, b=b$, a loop
labelled $1+T_i$ iff $(1+T_i)\, b=b$, and an edge from $b$ to $b'$ labelled
$T_i$ iff $T_i\, b_\tau$ is $b_{\s_i\tau}$
(if $\s_i\tau$ occurs in $\sigma\convol\sigma'$, and $\s_i\tau$ has one
inversion more than $\tau$) or $-b_\tau$ or $0$, depending on whether $i$ is a
descent of $\tau^{-1}$ or not.
The composition factors are then the simple modules indexed by the
descent compositions of the inverses of the permutations occuring in the
convolution~\cite{NCSF3, NCSF4} (see Figure~\ref{0-Hecke-graphe}).

This can be interpreted as a computation in the algebra
$\FQSym^*$~\cite{MR,NCSF6}, where the product of two basis elements
$\G_\alpha \G_\beta$ is
\begin{equation}
\G_\alpha \G_\beta = \sum_{\gamma\,\in\,\alpha\convol\beta} \G_\gamma.
\end{equation}
Indeed, the composition factors of the induction product of
two simple $0$-Hecke modules can be obtained by computing the product
$\G_\sigma \G_{\sigma'}$ and taking the image of the result by the morphism
sending $\G_\tau$ to the quasi-symmetric function
$F_{\Des(\tau^{-1})}\in\QSym$.
It is known that this amounts to take the commutative image of
the realization of $\FQSym^*$ by noncommutative polynomials~\cite{NCSF6}.
Since this is an epimorphism, the Grothendieck ring $\mathcal G$ of the
tower of the $0$-Hecke algebras is isomorphic to $\QSym$.
We shall later use a similar argument to identify the Grothendieck ring of
the $0$-Ariki-Koike-Shoji algebras.

\subsubsection{Indecomposable projective modules}

\begin{figure}[ht]
\PetitTableau
\ifdraft

\rotateleft{
\newdimen\vcadre\vcadre=0.2cm 
\newdimen\hcadre\hcadre=0.2cm 
\def\GrTeXBox#1{\vbox{\vskip\vcadre\hbox{\hskip\hcadre%
      $\Tableau{#1\\ }$%
   \hskip\hcadre}\vskip\vcadre}}
\def\GrTeXBoxb#1{\vbox{\vskip\vcadre\hbox{\hskip\hcadre%
   $#1$\hskip\hcadre}\vskip\vcadre}}

  \def\loopRight{\ar @`{[]+<1.4cm,1.1cm>,[]+<1.4cm,-1.1cm>}}
  \def\loopLeft{\ar @`{[]+<-1.4cm,1.1cm>,[]+<-1.4cm,-1.1cm>}}

$
\vcenter{\xymatrix@R=1.8cm@C=+0.6cm{
   & *{\GrTeXBox{2 \\ 1& 4 \\\ & 3}} \loopRight^{{}^{-T_1}_{-T_3}}
     \ar @{->}[d]^{T_2} & \\
   & *{\GrTeXBox{3 \\ 1& 4 \\\ & 2}} \loopRight^{-T_2}
   \ar @{->}[dl]_{T_1} \ar @{->}[dr]^{T_3} & \\
     *{\GrTeXBox{3 \\ 2& 4 \\\ & 1}} \loopLeft_{{}^{-T_1}_{-T_2}}
      \ar @{->}[dr]_{T_3}
 & & *{\GrTeXBox{4 \\ 1& 3 \\\ & 2}}\loopRight^{{}^{-T_2}_{-T_3}}
      \ar @{->}[dl]^{T_1} \\
   & *{\GrTeXBox{4 \\ 2& 3 \\\ & 1}} \loopRight^{{}^{-T_1}_{-T_3}}& \\
}}
$
}
\fi
\caption{\label{proj121}The projective module $P_{121}$ of $H_4(0)$.}
\end{figure}
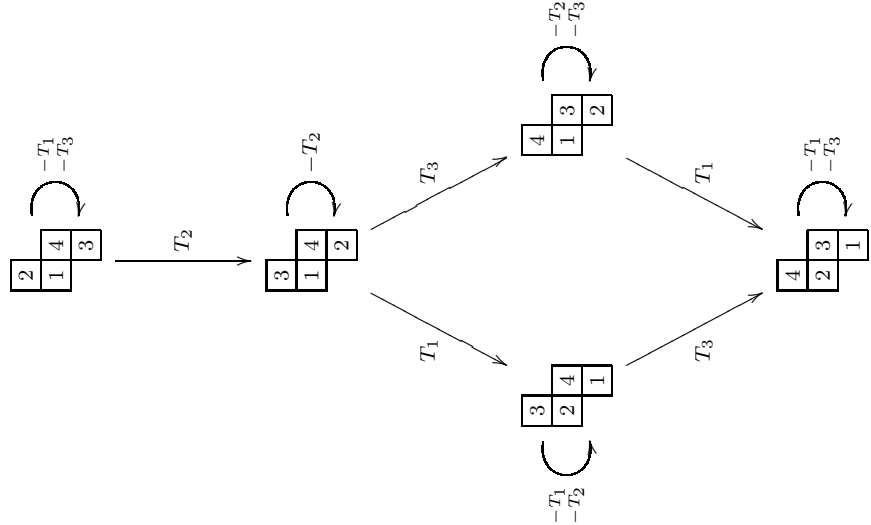

The indecomposable projective $H_n(0)$-mo\-dule $P_I$ labelled by a
composition $I$ of $n$ is combinatorial in a certain basis $(g_\sigma)$. Its
$0$-Hecke graph coincides with the restriction of the left weak order to the
interval $[\alpha(I),\omega(I)]$, that is the set of permutations with descent
composition $I$, an edge $\sigma\mapsto_i\tau$ meaning that
$T_i(g_\sigma)=g_\tau$.
The basis $g_\sigma$ is defined in terms of Norton's generators:
\begin{equation}
\label{Norton-nu}
g_{\alpha(I)} = \nu_I
= T_{\alpha(I)} T'_{\alpha(\overline{I}^\sim)},
\end{equation}
where $T'_i=1+T_i$ satisfies the braid relations, so that $T'_\sigma$ is well
defined for any permutation.
Figure~\ref{proj121} shows the projective module $P_{121}$ of $H_4(0)$.

\subsection{Yang-Baxter bases of $H_n(0)$}
\label{YB}

It is known that from a set of generators (depending on two parameters $t$ and
$u$)
\begin{equation}
Y_i(t,u) = a(t,u) + b(t,u) T_i\qquad \text{$b\not=0$},
\end{equation}
of $H_n(0)$ satisfying the quantum Yang-Baxter equation
\begin{equation}
\label{yb}
Y_i(t,u) Y_{i+1}(t,v) Y_i(u,v) = Y_{i+1}(u,v) Y_i(t,v) Y_{i+1}(t,u),
\end{equation}
one can associate with any vector $\x=(x_1,\ldots,x_n)$ of ``spectral
parameters'', a Yang-Baxter basis $(Y_\sigma(\x))_{\sigma\in\SG_n}$ of
$H_n(0)$, inductively defined by $Y_{id}(\x)=1$ and
\begin{equation}
Y_{\sigma_j\tau}(\x) = Y_j(x_{\tau^{-1}(j)},x_{\tau^{-1}(j+1)}) Y_{\tau}(\x)
= Y_{\sigma_j}(\x\tau^{-1}) Y_{\tau}(\x),
\end{equation}
if $\sigma_j\tau$ has one inversion more than $\tau$.
The special case $Y_j(t,u) = 1+ (1-u/t) T_j$ is studied in~\cite{LLT}.
Here we need the solution
\begin{equation}
\label{sgn}
Y_j(t,u) =
\left\{
\begin{array}{ll}
T_j   & \text{if $t  >  u$},\\
1+T_j & \text{if $t\leq u$},
\end{array}
\right.
\end{equation}
Indeed, checking the six different possibilities of order among $t$, $u$, $v$,
one has:
\begin{lemma}
\label{yb-lemme}
The $Y_j(t,u)$ defined by Equation~(\ref{sgn}) satisfy Equation~(\ref{yb}).
\qed
\end{lemma}

\begin{figure}[ht]
\ifdraft




\newdimen\vcadre\vcadre=0.20cm 
\newdimen\hcadre\hcadre=0.20cm 
\def\GrTeXBox#1{\vbox{\vskip\vcadre\hbox{\hskip\hcadre%
      $Y_{#1}$%
   \hskip\hcadre}\vskip\vcadre}}

\def\arx#1[#2](#3){\ifcase#1 \relax \or%
  \ar @{.>}[#2]|(#3)*+{\scriptstyle T_1}  \or%
  \ar @{->}[#2]|(#3)*+{\scriptstyle T_2} \or%
  \ar @{=>}[#2]|(#3)*+{\scriptstyle T_3} \or%
  \ar @{.>}[#2]|(#3)*+{\scriptstyle 1+T_1}  \or%
  \ar @{->}[#2]|(#3)*+{\scriptstyle 1+T_2} \or%
  \ar @{=>}[#2]|(#3)*+{\scriptstyle 1+T_3} \or%
\fi}

$\xymatrix@R=2.5cm@C=0.1cm{%
 &  &  &  &  & *{\GrTeXBox{1234}}\arx4[llld](.5)\arx2[d](.5)\arx3[rrrd](.5)& \\
 &  & *{\GrTeXBox{2134}}\arx5[ld](.5)\arx3[rrrd](.4)&
 &  & *{\GrTeXBox{1324}}\arx4[lld](.4)\arx3[rrd](.4)&
 &  & *{\GrTeXBox{1243}}\arx4[llld](.4)\arx2[rd](.5)& \\
 & *{\GrTeXBox{3124}}\arx1[rd](.5)\arx3[ld](.5)&
 & *{\GrTeXBox{2314}}\arx5[ld](.5)\arx3[rrrd](.4)&
 & *{\GrTeXBox{2143}}\arx2[ld](.35)&
 & *{\GrTeXBox{1423}}\arx4[ld](.5)\arx2[rrrd](.4)&
 & *{\GrTeXBox{1342}}\arx1[ld](.3)\arx3[rd](.5)& \\
   *{\GrTeXBox{4123}}\arx1[rd](.5)\arx2[rrrd](.35)&
 & *{\GrTeXBox{3214}}\arx3[ld](.35)&
 & *{\GrTeXBox{3142}}\arx1[rrrd](.35)\arx6[ld](.5)&
 & *{\GrTeXBox{2413}}\arx2[ld](.35)&
 & *{\GrTeXBox{2341}}\arx5[ld](.5)\arx3[rd](.5)&
 & *{\GrTeXBox{1432}}\arx1[ld](.5)& \\
 & *{\GrTeXBox{4213}}\arx2[rd](.5)&
 & *{\GrTeXBox{4132}}\arx1[rrd](.6)&
 & *{\GrTeXBox{3412}}\arx1[rrrd](.6)\arx6[llld](.6)&
 & *{\GrTeXBox{3241}}\arx6[lld](.6)&
 & *{\GrTeXBox{2431}}\arx5[ld](.5)& \\
 &  & *{\GrTeXBox{4312}}\arx1[rrrd](.5)&
 &  & *{\GrTeXBox{4231}}\arx2[d](.5)&
 &  & *{\GrTeXBox{3421}}\arx6[llld](.5)& \\
 &  &  &  &  & *{\GrTeXBox{4321}}& \\
}$
\fi
\caption{\label{permuto2431}The Yang-Baxter graph associated with the color
word $2431$.}
\end{figure}
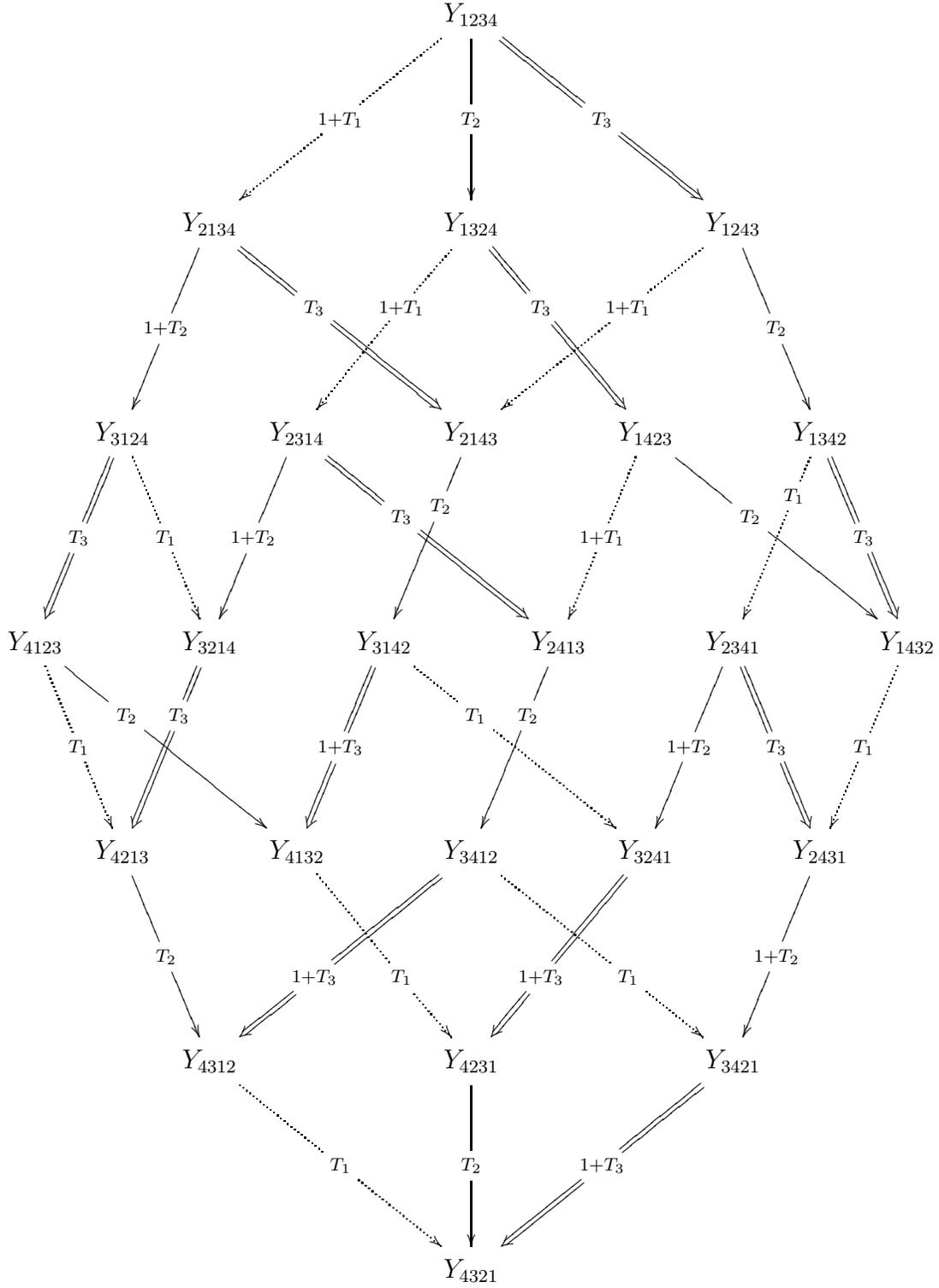

The \emph{Yang-Baxter graph} associated with a word $\x$ is the graph whose
vertices are $Y_\sigma(\x)$ and whose edges go from
$Y_\tau$ to $Y_{\sigma_j\tau}$ with the label
$Y_j(x_{\tau(j)},x_{\tau(j+1)})$ if $\sigma_j\tau$ has one inversion more than
$\tau$.
For example, one can check on Figure~\ref{permuto2431} that
$Y_{3241}(2431) = (1+T_2)T_1T_2T_3 = T_1T_2(1+T_1)T_3 = T_1T_2T_3(1+T_1)$.

\begin{note}
{\rm
The conditions of Equation~(\ref{sgn}) show that $Y_\sigma(\x)$ does not
depend on $\x$ but only on its standardized word $\Std(\x)$:
\begin{equation}
\label{Ystd}
Y_\sigma(\x) = Y_\sigma(\Std(\x)).
\end{equation}
We can therefore assume that $\x$ is a color word with no repeated colors,
hence a permutation.
}
\end{note}

\subsection{Yang-Baxter graphs of semi-combinatorial modules}
\label{propsYB}

We shall now generalize the definition of a combinatorial module as follows:
a \emph{semi-combinatorial $H_n(0)$-module} is a $H_n(0)$-module with a basis
$(b)$ satisfying: either $T_ib$, $-T_ib$ or $(1+T_i)b$ is a basis element.

\begin{figure}[ht]
\PetitTableau

\rotateleft{
\newdimen\vcadre\vcadre=0.2cm 
\newdimen\hcadre\hcadre=0.0cm 
\def\GrTeXBox#1{\vbox{\vskip\vcadre\hbox{\hskip\hcadre%
      $\Tableau{#1\\ }$%
   \hskip\hcadre}\vskip\vcadre}}
\def\GrTeXBoxb#1{\vbox{\vskip\vcadre\hbox{\hskip\hcadre%
   $#1$\hskip\hcadre}\vskip\vcadre}}
  \def\loopRight{\ar @`{[]+<1.2cm,1.1cm>,[]+<1.2cm,-1.1cm>}}
  \def\loopLeft{\ar @`{[]+<-1.2cm,1.1cm>,[]+<-1.2cm,-1.1cm>}}
  \def\loopRightb{\ar @`{[]+<1.5cm,1.1cm>,[]+<1.5cm,-1.1cm>}}
  \def\loopLeftb{\ar @`{[]+<-1.8cm,1.1cm>,[]+<-1.8cm,-1.1cm>}}

$\vcenter{\xymatrix@R=1.6cm@C=+0.3cm{
   & *{\GrTeXBoxb{Y_{2143}}} \ar @{->}[d]^{1+T_2} & \\
   & *{\GrTeXBoxb{Y_{3142}}}
   \ar @{..>}[dl]_{1+T_1} \ar @{=>}[dr]^{T_3}& \\ 
 *{\GrTeXBoxb{Y_{3241}}}\ar @{=>}[dr]_{T_3}
 & & *{\GrTeXBoxb{Y_{4132}}} \ar @{..>}[dl]^{1+T_1}\\
   & *{\GrTeXBoxb{Y_{4231}}} \ar @{->}[d]^{T_2}& \\
   & *{\GrTeXBoxb{Y_{4321}}} & \\
}}
\quad\longmapsto\quad
\vcenter{\xymatrix@R=1.6cm@C=+0.3cm{
   & *{\GrTeXBoxb{Y_{2143}}} \loopRight^{{}^{1+T_1}_{1+T_3}}  \ar @{->}[d]^{1+T_2} & \\
   & *{\GrTeXBoxb{Y_{3142}}} \loopRight^{{}^{1+T_2}}
   \ar @{..>}[dl]_{1+T_1} \ar @{=>}[dr]^{T_3}& \\ 
 *{\GrTeXBoxb{Y_{3241}}}\loopLeft_{{}^{1+T_1}_{1+T_2}} \ar @{=>}[dr]_{T_3}
 & & *{\GrTeXBoxb{Y_{4132}}}\loopRight^{{}^{1+T_2}_{-T_3}} \ar @{..>}[dl]^{1+T_1}\\
   & *{\GrTeXBoxb{Y_{4231}}}\loopRight^{{}^{1+T_1}_{-T_3}} \ar @{->}[d]^{T_2}& \\
   & *{\GrTeXBoxb{Y_{4321}}}\loopRight^{{}^{1+T_1}_{-T_2,1+T_3}} & \\
}}
\quad\longmapsto\quad
\vcenter{\xymatrix@R=1.2cm@C=+0.2cm{
   & *{\GrTeXBox{ &  \\\ &  & }} \ar @{->}[d] & \\
   & *{\GrTeXBox{ \\ & & }}
   \ar @{->}[dl] \ar @{->}[dr]& \\
 *{\GrTeXBox{&  &  & }} \ar @{->}[dr]
 & & *{\GrTeXBox{\\ & \\\ &}} \ar @{->}[dl]\\
   & *{\GrTeXBox{&& \\\ &\ &}} \ar @{->}[d]&
\\
   & *{\GrTeXBox{& \\\ & &}} & \\
}}
$
}
\caption{\label{fig-yb}The first graph is the restriction of the Yang-Baxter
graph associated with the color word $2314$ to the convolution $21\convol21$.
The second graph describes the corresponding semi-combinatorial
$H_n(0)$-module. In the third graph, each basis element $Y_\sigma$ has been
replaced by the associated composition factor, given by the shape of the
anticycloribbon associated with $(\sigma,2314)^{-1}$.}
\end{figure}
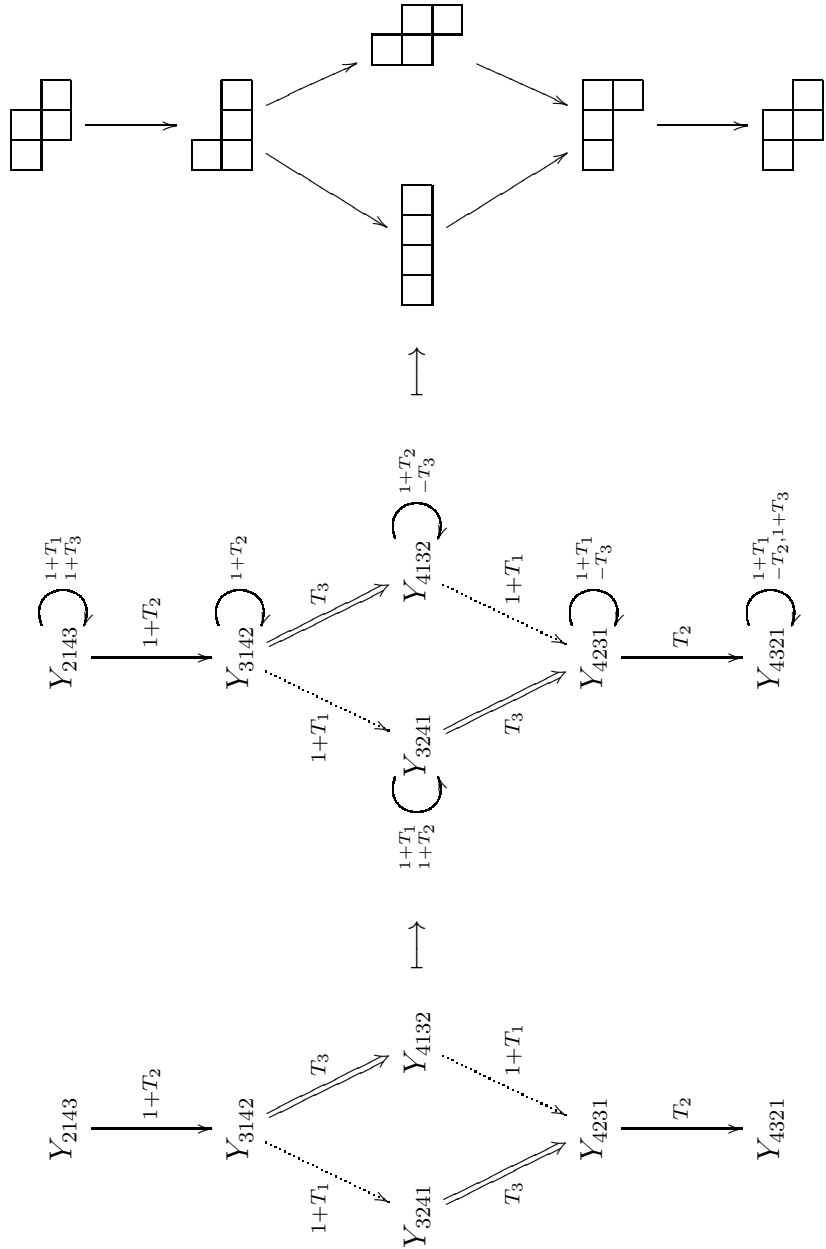

In particular, from any interval $[Y_\sigma(\x), Y_\tau(\x)]$ in a
Yang-Baxter graph, one can build the graph of a semi-combinatorial module of
$H_n(0)$, in the basis $(Y_\rho)$, by adding loops on each vertex $Y_\rho$
as follows: if $Y_{\s_i\rho}$ does not belong to the interval, add a loop
on $Y_\rho$ labelled $-T_i$ or $1+T_i$ depending on whether
$\x_{\rho^{-1}(i)}<\x_{\rho^{-1}(i+1)}$ or not.
As in the case of combinatorial $H_n(0)$-modules, one can read on the graph
$G$ of a semi-combinatorial module as many composition series of the module as
linear extensions of $G$.
One can see an example of a restricted Yang-Baxter graph and the corresponding
semi-combinatorial module on Figure~\ref{fig-yb}.

Note that the interval $[Y_\sigma(\x), Y_{\omega_n}(\x)]$
corresponds to the left $H_n(0)$-module generated by $Y_\sigma(\x)$.
Consequently, any interval $[Y_\sigma(\x), Y_\tau(\x)]$ corresponds to a
quotient of the previous $H_n(0)$-module, hence to an $H_n(0)$-module itself.

\begin{lemma}
\label{simpleYB}
Let $\x$ be a word, $\sigma$ a permutation and consider the interval of the
left weak order $[Y_\sigma(\x),Y_\sigma(\x)]$ as a one-dimensional
quotient of the $H_n(0)$-module generated by $Y_\sigma(\x)$.
It is isomorphic to the simple $H_n(0)$-module labelled by the shape of the
anticycloribbon associated with $(\sigma,\x)^{-1}$.

In particular, $Y_{\omega_n}(\x)$ is the simple $H_n(0)$-module labelled by
the shape of the anticycloribbon associated with $(\omega_n,\overline{\x})$.
\end{lemma}

\begin{proof}
The action of $T_i$ on the quotient $[Y_\sigma(\x),Y_\sigma(\x)]$ is
encoded by the Yang-Baxter graph associated with $\x$, and the
lemma follows from the definitions.
\end{proof}

For example, the one-dimensional quotient of the left $H_n(0)$-module
generated by $Y_{1342}(3213)$ has as corresponding anticycloribbon
\begin{equation}
\PetitTableau
\Tableau{1^3\\ 4^3 & 2^2 & 3^1\\ }
\end{equation}
so that it is the simple module $S_{13}$, as one can check on
Figure~\ref{fig-yb} (the generator of the quotient module is annihilated by
$1+T_1$, $-T_2$ and $-T_3$).

\subsection{Representations in Yang-Baxter bases}

\subsubsection{Simple modules and indecomposable projective modules}

We can now rewrite the representation theory of $H_n(0)$ in terms of
Yang-Baxter bases. First, we identify the simple and indecomposable projective
modules of $H_n(0)$ as the semi-combinatorial modules associated with
some intervals of Yang-Baxter graphs. Then, using Theorem~\ref{induit-shuf},
we will be able to describe in a very simple way the induction of those
modules in terms of intervals in the Yang-Baxter graphs.

We first consider the case of simple $H_n(0)$-modules.
The next lemma is a rewriting of a special case of Lemma~\ref{simpleYB}.
Recall that $\eta_I$ is the generator of the simple module $S_I$ (see
Equation~(\ref{etaI})).

\begin{lemma}
\label{omega-eta}
Let $I$ be a composition. Then for any permutation $\tau$ such that
$\overline{\Des(\tau)}=I$, one has
\begin{equation}
\eta_I = Y_{\omega_n}(\tau).
\end{equation}
\end{lemma}

Let $\antiinvol$ be the anti-automorphism of $H_n(0)$ such that
$\antiinvol(T_i)=T_i$ for all $i$.
One easily proves by induction on the number of inversions of $\sigma$ that

\begin{lemma}
\label{ech-hecke}
For all color words $\tau$ and all permutations $\sigma$, one has
\begin{equation}
\antiinvol(Y_\sigma(\tau))
= Y_{\sigma^{-1}}(\overline{\tau}^\#\cdot(\sigma^{-1})^\#).
\end{equation}
\end{lemma}

These lemmas give in particular the description of $\CC \eta_I$
as a right $H_n(0)$-module:
\begin{proposition}
\label{ech-eta}
For all color words $\tau$, one has
\begin{equation}
\antiinvol(Y_{\omega_n}(\tau)) = Y_{\omega_n}(\overline{\tau}^\#)),
\end{equation}
so that, for all compositions $I$,
\begin{equation}
\antiinvol(\eta_I) = \eta_{\overline{I}}.
\end{equation}
\end{proposition}

Let us now consider the case of indecomposable projective $H_n(0)$-modules.
For those modules, from the definition of Norton's generators (see
Equation~(\ref{Norton-nu})), one obtains:

\begin{proposition}
Let $I$ be a composition. Then
%
\begin{equation}
\nu_I = Y_{\alpha(I)\cdot\alpha(I^\sim)}(\omega(I)).
\end{equation}
\end{proposition}

\subsubsection{Induction product of modules}

\begin{theorem}
\label{induit-shuf}
Let $M'$ and $M''$ be the semi-combinatorial modules of $H_k(0)$ and
$H_{n-k}(0)$ associated with the intervals
$[Y_{\alpha'}(\tau'),Y_{\beta'}(\tau')]$ and
$[Y_{\alpha''}(\tau''),Y_{\beta''}(\tau'')]$.
Then, $M'\widehat\otimes M''$ is the semi-combinatorial module associated with
the interval
$[Y_{\alpha}(\tau),Y_{\beta}(\tau)]$, where
$\tau$ is any element of $\tau'\convol\tau''$,
$\alpha=\alpha'\cdot\alpha''[k]$ and $\beta=\beta''[k]\cdot\beta'$.
\end{theorem}

\begin{proof}
First, the semi-combinatorial module $M$ corresponding to the restriction of
the Yang-Baxter graph determined by $\tau$ to the interval
$[\alpha,\beta]$ and $M'\widehat\otimes M''$ have the same dimension and both
are $H_n(0)$-modules.
Moreover, the quotient of $M$ by the submodule $H_n(0) T_k M$ is isomorphic
to $M'\otimes M''$ as an $H_k(0)\otimes H_{n-k}(0)$-module. So $M$ is
isomorphic to $M'\widehat\otimes M''$.
\end{proof}

Since any simple or indecomposable projective module is an interval of a
Yang-Baxter graph (seen as a semi-combinatorial module), the induction of such
modules is also such a Yang-Baxter graph. The choice we have in color $\tau$
in Theorem~\ref{induit-shuf} will be useful when computing the induced module
of two indecomposable projective modules in the $0$-Ariki-Koike-Shoji
algebras.

\begin{corollary}
\label{yb-bases}
Let $I_1$ and $I_2$ be compositions of $n_1$ and $n_2$, and
$\x=\x_1\concat\x_2$ a color word with
$\x_1\in C^{n_1}$ and $\x_2\in C^{n_2}$.
Let $\sigma_1=\omega(I_1)$ and $\sigma_2=\omega(I_2)$.

Then, the semi-combinatorial module corresponding to the restriction of the
Yang-Baxter graph determined by $\x$ to the convolution
$\sigma_1\convol\sigma_2$ is isomorphic to the induced module
$S_{I'_1}\widehat\otimes S_{I'_2}$ where $I'_1$ (resp. $I'_2$) is the shape of
the anticycloribbon associated with the colored permutation
$(\sigma_1,\x_1)^{-1}$
(resp. $(\sigma_2,\x_2)^{-1})$.
\end{corollary}

\begin{corollary}
\label{labels-induc-simples}
Let $I'$ and $I''$ be two compositions of integers $n'$ and $n''$. Let $\x'$
and $\x''$ be two permutations such that $I'=\overline{C(\x')}$ and
$I''=\overline{C(\x'')}$.
Then the interval $[Y_\sigma^\x,Y_{\omega_{n'+n''}}^\x]$ of the Yang-Baxter
graph of color $\x$ where $\sigma=\omega_{n'}\sconc\omega_{n''}$ and
$\x=\x'\sconc\x''$ is isomorphic the induced module
$S_{I'}\widehat\otimes S_{I''}$.
\end{corollary}

For example, the semi-combinatorial module represented on Figure~\ref{fig-yb} 
is isomorphic to $S_{2}\widehat{\otimes} S_{2}$, as one can
check on Figure~\ref{0-Hecke-graphe}. The edges are not labelled by the same
elements but the composition factors are the same. Since the generators of
both modules are the same, one can easily provide an explicit description of
the basis of the semi-combinatorial module in the left regular representation:
it is a product of factors $T_i$ and $1+T_i$ multiplied on their right by
$\eta_{I'_1}\otimes\eta_{I'_2}$ seen as an element of $H_{|I_1|+|I_2|}(0)$
through the morphism $H_{|I_1|}(0)\otimes H_{|I_2|}(0)\to H_{|I_1|+|I_2|}(0)$.

\begin{figure}[ht]
\PetitTableau

{
\newdimen\vcadre\vcadre=0.2cm 
\newdimen\hcadre\hcadre=0.0cm 
\def\GrTeXBox#1{\vbox{\vskip\vcadre\hbox{\hskip\hcadre%
      $\Tableau{#1\\}$%
   \hskip\hcadre}\vskip\vcadre}}
\def\GrTeXBoxb#1{\vbox{\vskip\vcadre\hbox{\hskip\hcadre%
   $#1$\hskip\hcadre}\vskip\vcadre}}
  \def\loopRight{\ar @`{[]+<1.2cm,1.1cm>,[]+<1.2cm,-1.1cm>}}
  \def\loopLeft{\ar @`{[]+<-1.2cm,1.1cm>,[]+<-1.2cm,-1.1cm>}}
  \def\loopRightb{\ar @`{[]+<1.5cm,1.1cm>,[]+<1.5cm,-1.1cm>}}
  \def\loopLeftb{\ar @`{[]+<-1.8cm,1.1cm>,[]+<-1.8cm,-1.1cm>}}

$
\vcenter{\xymatrix@R=1.0cm@C=+0.3cm{
   & *{\GrTeXBoxb{Y_{2314}}} \loopRight^{{}^{-T_1}}  \ar @{->}[dl]_{T_2} \ar @{=>}[dr]^{1+T_3}& \\
 *{\GrTeXBoxb{Y_{3214}}}\loopLeft_{{}^{1+T_1}_{-T_2}} \ar @{=>}[d]_{1+T_3}
 & & *{\GrTeXBoxb{Y_{2413}}}\loopRight^{{}^{-T_1}_{1+T_3}} \ar @{->}[d]^{1+T_2}\\
 *{\GrTeXBoxb{Y_{4213}}}\loopLeft_{{}^{1+T_1}_{1+T_3}} \ar @{->}[dr]_{1+T_2}
 & & *{\GrTeXBoxb{Y_{3412}}}\loopRight^{{}^{1+T_2}} \ar
 @{=>}[dl]^{T_3} \ar @{..>}[dr]^{1+T_1}\\
   & *{\GrTeXBoxb{Y_{4312}}}\loopLeft_{{}^{1+T_2}_{-T_3}} \ar @{..>}[dr]_{1+T_1}
   & & *{\GrTeXBoxb{Y_{3421}}}\loopRight^{{}^{1+T_1}_{-T_2}} \ar @{=>}[dl]^{T_3} \\
 & & *{\GrTeXBoxb{Y_{4321}}}\loopRight^{{}^{1+T_1}_{1+T_2,-T_3}} & \\
}}
\hfill
\vcenter{\xymatrix@R=1.0cm@C=+0.3cm{
   & *{\GrTeXBoxb{Y_{2314}}} \loopRight^{{}^{-T_1}}  \ar @{->}[dl]_{T_2} \ar @{=>}[dr]^{T_3}& \\
 *{\GrTeXBoxb{Y_{3214}}}\loopLeft_{{}^{1+T_1}_{-T_2}} \ar @{=>}[d]_{T_3}
 & & *{\GrTeXBoxb{Y_{2413}}}\loopRight^{{}^{-T_1}_{-T_3}} \ar @{->}[d]^{T_2}\\
 *{\GrTeXBoxb{Y_{4213}}}\loopLeft_{{}^{1+T_1}_{-T_3}} \ar @{->}[dr]^{T_2}
 & & *{\GrTeXBoxb{Y_{3412}}}\loopRight^{{}^{-T_2}} \ar
 @{=>}[dl]_{T_3} \ar @{..>}[dr]^{T_1}\\
   & *{\GrTeXBoxb{Y_{4312}}}\loopLeft_{{}^{-T_2}_{-T_3}} \ar @{..>}[dr]_{T_1}
   & & *{\GrTeXBoxb{Y_{3421}}}\loopRight^{{}^{-T_1}_{-T_2}} \ar @{=>}[dl]^{T_3} \\
 & & *{\GrTeXBoxb{Y_{4321}}}\loopRight^{{}^{-T_1}_{1+T_2,-T_3}} & \\
}}
$
}
\caption{\label{fig-projind}
The graphs are the semi-combinatorial modules corresponding to the restriction
of the Yang-Baxter graph associated with the color words $3124$ and $4231$ to
the same interval $[2314,4321]$. Both are the graph of the induced module 
$P_{(2,1)}\widehat\otimes P_{(1)}$}.
\end{figure}
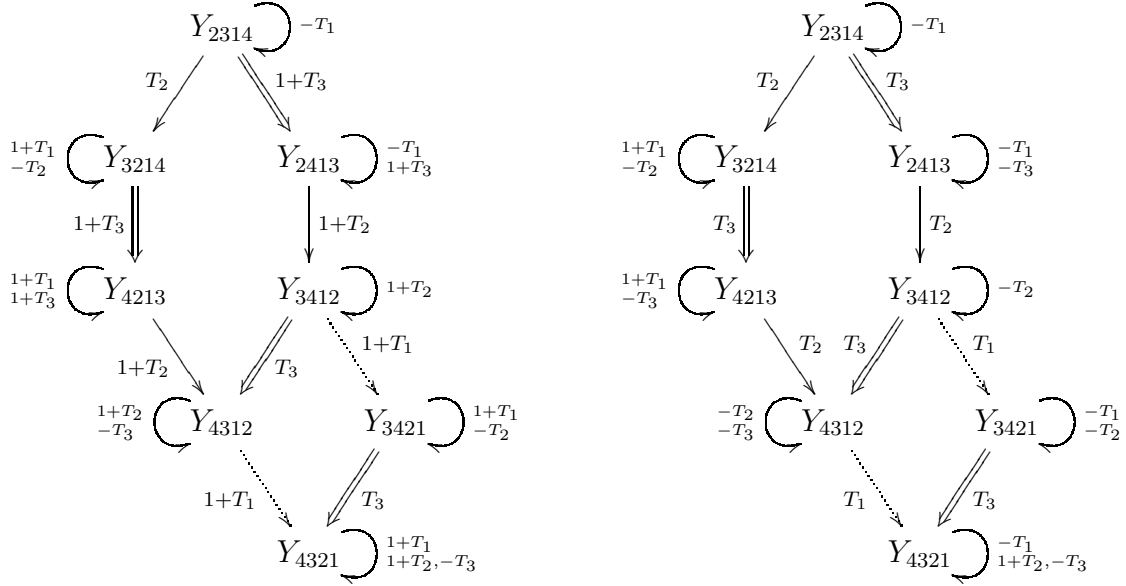

Another example is given Figure~\ref{fig-projind}. The graphs are the
semi-combinatorial modules corresponding to the restriction of 
the Yang-Baxter graph associated with the color words $3124$ and $4231$ to the
same interval $[2314,4321]$. Both are the graph of the induced module
$P_{(2,1)}\widehat\otimes P_{(1)}$ where $P_I$ is the indecomposable
projective module associated with $I$.

\newpage
\section{The $0$-Ariki-Koike-Shoji algebras}

In \cite{Sho}, Shoji obtained a new presentation of the Ariki-Koike algebras
defined in~\cite{AK}.
We shall first give a presentation very close to his, put $q=0$ in the
relations and then prove some simple results about another basis of the
resulting algebra. To get our presentation from Shoji's~\cite{Sho}, one has to
replace $qa_i$ by $T_{i-1}$ and $q^2$ by $q$.

Let $u_1,\ldots,u_r$ be $r$ distinct complex numbers. We shall denote by
$P_k(X)$ the Lagrange polynomial
\begin{equation}
\label{Lag}
P_k(X) := \prod_{1\leq l\leq r,\ l\not=k} \frac{X-u_l}{u_k-u_l}.
\end{equation}

The Ariki-Koike algebra $\AKSq$ is the associative $\CC$-algebra generated by
elements $T_1, \ldots,T_{n-1}$ and $\xi_1,\ldots,\xi_n$ subject to the
following relations:

\begin{alignat}{2}
(T_i-q)(T_i+1)        &= 0                    && (1\leq i\leq n-1), \\
T_i T_{i+1} T_i       &= T_{i+1} T_i T_{i+1}  && (1\leq i\leq n-2), \\
T_i T_j               &= T_j T_i              && (|i-j|\geq2), \\
{(\xi_j-u_1)\cdots (\xi_j-u_r)} &= 0          && (1\leq j\leq n), \\
\xi_i \xi_j           &= \xi_j \xi_i  && (1\leq i,j \leq n), \label{comxi} \\
T_i \xi_i = \xi_{i+1}T_i - (q-1)
&\sum_{c_1<c_2}{(u_{c_2}-u_{c_1})P_{c_1}(\xi_i)P_{c_2}(\xi_{i+1})}
&\quad& (1\leq i\leq n-1), \label{tx-xt}\\
T_i (\xi_{i+1}+\xi_i) &= (\xi_{i+1}+\xi_i)T_i && (1\leq i\leq n-1)
\label{txx-xxt}\\
T_i \xi_j             &= \xi_j T_i            && (j\not= i+1,i).
\end{alignat}

As noticed in~\cite{Sho}, it is obvious from this presentation that a
generating set is given by the $\xi_1^{c_1}\cdots\xi_n^{c_n}\cdot T_\sigma$
with $\sigma\in\SG_n$ and $c_i$ such that $0\leq c_i\leq r-1$. Shoji proves
that this is indeed a basis of $\AKSq$.
A simple adaptation of his proof would enable us to conclude that this
property still holds at $q=0$.
We will prove it by introducing a new basis on which the product
of generators has a simple expression (see Lemma~\ref{commut}.)
This algebra $\AKSz$, which we call the $0$-Ariki-Koike-Shoji algebra, will be
our main concern in the sequel.

If $\c=(c_1,\ldots,c_n)$ is a word on $C=\{1,\ldots,r\}$, we define
\begin{equation}
L_{\c} := P_{c_1}(\xi_1) \cdots P_{c_n}(\xi_n).
\end{equation}

Since the Lagrange polynomials (\ref{Lag}) associated with $r$
distincts complex numbers form a basis of $\CC_{r-1}[X]$ (polynomials of
degree at most $r-1$), the set of elements of $\AKSq$
\begin{equation}
\label{baseB}
B := \{ B_{\c,\sigma}(q) := L_{\c} T_\sigma\},
\end{equation}
where $\sigma$ runs over $\SG_n$ and $\c=(c_1,\ldots,c_n)$ runs over
the color words of size $n$ span the same vector space as the
$\xi^{\c}T_\sigma$.

Let us now describe the action of a generator on the left of
$B_{\c,\sigma}$: the generator $\xi_i$ acts diagonally by
multiplication by $u_{c_i}$, so that it only remains to explicit the product
$T_i \, L_{\c} \, T_\sigma$.
The relevant expression comes from the following apparently unnoticed relation
in
$\AKSq$:
\begin{lemma}
\label{commut}
The following relation holds in $\AKSq$:
\begin{equation}
\label{commutQ}
T_i L_{\c} = L_{\c\s_i} T_i - (q-1)\left\{
\begin{array}{lcl}
-L_{\c} & \text{if} &c_i<c_{i+1}, \\
0 & \text{if} &c_i=c_{i+1}, \\
L_{\c\s_i} & \text{if} &c_i>c_{i+1}. \\
\end{array}
\right.
\end{equation}
where $\s_i$ acts on the right of $\c$ by exchanging $c_i$ and $c_{i+1}$.
\end{lemma}

\begin{proof}
It is enough to compute $T_iP_{a}(\xi_i)P_{b}(\xi_{i+1})$ since all the other
terms commute with $T_i$. Using Relations~(\ref{comxi}),~(\ref{tx-xt}),
and~(\ref{txx-xxt}), one proves
\begin{equation}
T_i\xi_i\xi_{i+1} = \xi_i\xi_{i+1}T_i,
\end{equation}
so that $T_i$ commutes with any symmetric function of $\xi_i$ and $\xi_{i+1}$.
The computation then reduces to evaluate the expression
\begin{equation}
\label{calculcommut}
\prod_{1\leq l\leq r,\ l\not=a,\ l\not=b}
\left({\frac{\xi_i-u_l}{u_a-u_l}}\ {\frac{\xi_{i+1}-u_l}{u_a-u_l}} \right)
\ \ \
T_i\
\frac{\xi_i-u_b}{u_a-u_b}
\frac{\xi_{i+1}-u_a}{u_b-u_a}.
\end{equation}
Now, writing $(\xi_i-u_b)(\xi_{i+1}-u_a)$ as
$\xi_i\xi_{i+1} - u_b(\xi_i+\xi_{i+1}) + u_au_b + (u_b-u_a)\xi_i$,
the first three terms commute with $T_i$ and the last one, thanks to
Equation~(\ref{tx-xt}), simplifies into a sum of two terms since all but one
Lagrange polynomial vanish when multiplied with the left factor of
Equation~(\ref{calculcommut}), the terms depending on the order relation
between $a$ and $b$.
\end{proof}

The last lemma implies that the set $B$ defined in Equation~(\ref{baseB}) is a
generating set of $\AKSq$. We now show that it is a basis of $\AKSq$ for
all values of $q$ and in particular for $q=0$.
One can argue as in~\cite{Sho}.
Let $V=\bigoplus_{i=1}^r{V_i}$, with $V_i=\CC^{m_i}$, $m_i\geq n$, and let
$W = \{v_1,\ldots,v_m\}$ with $m=m_1+\cdots+m_r$, be a basis of $V$ such that
$v_1,\ldots,v_{m_1}\in V_1$, $v_{m_1+1},\ldots,v_{m_1+m_2}\in V_2$, and so on.
Let $w=v_{k_1}\otimes\cdots\otimes v_{k_n}\in V^{\otimes n}$.
There is a classical right action of $H_n(q)$ on $V^{\otimes n}$ given by

\begin{equation}
\label{eq:hecktensor}
\left\{
\begin{array}{cccl}
w \cdot T_i & = & w\s_i & \qquad
\hbox{if} \ k_i < k_{i+1} \ , \\
w \cdot T_i & = & q\, w & \qquad
\hbox{if} \ k_i = k_{i+1} \ , \\
w \cdot T_i & = & q\,w {\s_i}
+ (q-1)w & \qquad \hbox{if} \ k_i > k_{i+1} \ . \\
\end{array}
\right.
\end{equation}
Let $\xi$ be the linear map whose restriction to $V_i$ is $u_i\cdot id_{v_i}$,
and $\xi_j=id^{\otimes(j-1)} \otimes \xi \otimes id^{\otimes(n-j)}$ acting on
the right of $V^{\otimes n}$.
Then, as observed by Shoji, the $T_i$ and $\xi_j$ generate a right action of
$\AKSq$ on $V^{\otimes n}$.
It is known that the representation~(\ref{eq:hecktensor}) of $H_n(q)$ is
faithful for all $q$, provided that $m_i\geq n$. This easily implies the next
proposition:

\begin{proposition}
The set $B$ defined in Equation~(\ref{baseB}) is a basis of $\AKSq$, for all
values of $q$.
\end{proposition}

With exactly the same arguments, one would also prove:
\begin{proposition}
\label{ensBp}
The set $B'=\{B'_{\c,\sigma}:=T_\sigma L_{\c}\}$ is a basis of
$\AKSq$, for all values of $q$.
\end{proposition}

We now provide a new presentation of $\AKSz$. Since $\AKSz$ satisfies the
relations of the next proposition, a simple comparison of dimensions proves
that:

\begin{proposition}
\label{notrepres}
The algebra with generators $L_{\c}$, $\c$ being any color word of size $r$,
and $T_i$ with $1\leq i\leq n-1$, and relations
\begin{eqnarray}
T_i(1+T_i) = 0 \quad (1\leq i\leq n-1), \label{ti1}\\
T_i T_{i+1} T_i = T_{i+1} T_i T_{i+1} \quad (1\leq i\leq n-2), \label{ti2}\\
T_i T_j = T_j T_i \qquad (|i-j|\geq2),\label{ti3} \\
L_{\bf d} L_{\c} = \delta_{{\bf d},\c} L_{\c},\hskip3cm
\label{LL}\\
\left\{
\begin{array}{rcl}
\label{TL}
\text{if} &c_i<c_{i+1},& (1+T_i) L_{\c} = L_{\c\s_i} T_i, \\[3pt]
\text{if} &c_i=c_{i+1},& T_i L_{\c} = L_{\c\s_i} T_i, \\[3pt]
\text{if} &c_i>c_{i+1},& T_i L_{\c} = L_{\c\s_i} (1+T_{i}). \\
\end{array}
\right.
\end{eqnarray}
is isomorphic to $\AKSz$.
\end{proposition}

This description of the Ariki-Koike-Shoji algebra enables us to analyze the
left regular representation of $\AKSz$ in terms of our basis elements. As a
first application, we shall obtain a classification of the simple
$\AKSz$-modules.

\subsection{The associative bilinear form}

Recall that a bilinear form $(\,,\,)$ on an algebra $A$ is said to
be associative if $(ab,c)=(a,bc)$ for all $a,b,c\in A$, and that $A$
is called a {\em Frobenius algebra} whenever it has a nondegenerate
associative bilinear form. Such a form induces an isomorphism of left
$A$-modules between $A$ and the dual $A^*$ of the right regular
representation. Frobenius algebras are in particular self-injective,
so that finitely generated projective and injective modules coincide
(see \cite{CR}).

It is known that $H_n(0)$ is a Frobenius algebra for the following bilinear
form:
for a basis $(Y_\sigma)$, we denote by $(Y^*_\sigma)$ the dual basis. We set
$\chi=T^*_{\omega_n}$. Then $(f,g)=\chi(fg)$ is a non-degenerate associative
bilinear form \cite{NCSF6}. Moreover if we denote $\eta_i=-T_i$ and
$\xi_i = 1+T_i$ (both satisfy the braid relations) then the following
properties hold:
\begin{equation}
\zeta_\sigma=(-1)^{\ell(\omega\sigma^{-1})}\xi_{\omega\sigma^{-1}},
\quad
\xi_\alpha =\sum_{\beta\le\alpha}T_\beta    \quad \text{and} \quad
\eta_{\alpha}=\sum_{\beta\le\alpha}(-1)^{\ell(\beta)}\xi_\beta
\end{equation}
shows that $(\zeta_\sigma)$ and $(\eta_\tau)$ are two adjoint bases of
$H_n(0)$, that is
\begin{equation}
(\zeta_\sigma,\eta_\tau)=\delta_{\sigma,\tau}\,.
\end{equation}

This bilinear form can be extended to $\AKSz$ by setting
\begin{equation}
\Chi := \sum_{\c} B_{\c,\omega}^*
\end{equation}
where $B$ is the basis defined by Equation~(\ref{baseB}).

\begin{proposition}
The associative bilinear form defined by
\begin{equation}
(f,g):=\Chi(fg)
\end{equation}
is non-degenerate on $\AKSz$. Therefore, $\AKSz$ is a Frobenius algebra.
\end{proposition}

\begin{proof}
The form is obviously bilinear and associative. It remains to prove that it
is non-degenerate. We have 
\begin{equation}
L_\c T_\sigma = T_\sigma L_{\c\sigma^{-1}} + \text{smaller terms},
\end{equation}
so that
\begin{equation}
(L_{\c}T_\sigma, T_\tau L_\d) = \Chi(L_{\c}T_\sigma T_\tau L_\d) =
\delta_{\c \overline\d} \Chi(T_\sigma T_\tau)\,.
\end{equation}
This formula being linear on $T_\sigma$ and $T_\tau$, we have as well,
\begin{equation}
\Chi(L_{\c}\xi_\alpha \eta_\beta L_{\overline\d}) =
\delta_{\c \d} \delta_{\alpha \beta}\,,
\end{equation}
so that $B'_{\c,\alpha} := L_\c \xi_\alpha$ and $B''_{\d,\beta} := \eta_\beta
L_{\overline d}$ are two adjoint bases. 
\end{proof}

In particular, we have
\begin{corollary}
\label{autoinj}
$\AKSz$ is self-injective.
\qed
\end{corollary}

\newpage
\section{Simple modules of $\AKSz$}

Let $[I,\c]$ be a cycloribbon.
Let $\eta_{[I,\c]}$ be the element of $\AKSz$ defined as
\begin{equation}
\eta_{[I,\c]} := L_{\c}\ \eta_I,
\end{equation}
where $\eta_I$ is the generator of the simple $H_n(0)$-module associated with
$I$ as in Section~\ref{rthn}.
This element generates a simple $\AKSz$-module:

\begin{theorem}
\label{simples}
Let $[I,\c]$ be a cycloribbon. Then
\begin{equation}
S_{[I,\c]} := \AKSz\, \eta_{[I,\c]}
\end{equation}
is a simple module of $\AKSz$, realized as a minimal left ideal in its left
regular representation. The eigenvalue of $L_{\d}$ is $1$ if $\d=\c$ and $0$
otherwise, and that of $T_i$ is $-1$ or $0$ according to whether $i$ is a
descent of the shape of $\bijcyc(R)$ or not.
All these simple modules are pairwise non isomorphic and of dimension $1$.
Moreover, all one-dimensional $\AKSz$-modules are isomorphic to
some $S_{[I,\c]}$.
\end{theorem}

\begin{proof}
We prove simultaneously that the $S_{[I,\c]}$ are one-dimensional
$\AKSz$-modules and that all one-dimensional modules are isomorphic to some
$S_{[I,\c]}$.
Let us first consider a simple module of dimension $1$. Since $\AKSz$ is
self-injective (Corollary~\ref{autoinj}), it can be realized as a minimal
left ideal in the left regular representation, that is, as $\AKSz S$ with
$S = \sum_{\c,\sigma} C_{\c,\sigma} B_{\c,\sigma}$ with $C_{\c,\sigma}\in\CC$.
Since all $B_{\c,\sigma}$ are eigenvectors of $L_\d$ with eigenvalue $1$ or
$0$ depending on whether $\d=\c$ or not, $S$ is an eigenvector for all $L_\d$
iff the only non-zero coefficients correspond to the same $\c$.
So $S = L_{\c} S'$ with $S'=\sum_{\sigma} C_{\sigma} T_\sigma$.
Since $S$ is also an eigenvector for all $T_i$, we are left with the following
possibilities:

\begin{itemize}
\item If $c_i<c_{i+1}$, then thanks to Formula~(\ref{TL}), first case,
$S$ is an eigenvector of $T_i$ iff $T_i$ acts by $0$ on $S'$.
\item If $c_i=c_{i+1}$, then thanks to Formula~(\ref{TL}), second case,
$S$ is an eigenvector of $T_i$ iff $T_i$ acts by $0$ or $-1$ on $S'$.
\item If $c_i>c_{i+1}$, then thanks to Formula~(\ref{TL}), third case,
$S$ is an eigenvector of $T_i$ iff $T_i$ acts by $-1$ on $S'$.
\end{itemize}
So, in any case, $S'$ is an eigenvector for all $T_i$, so it is equal
to some $\eta_I$. Moreover, $I$ has to be compatible with
$\c$ in the following sense: if $c_i<c_{i+1}$ then $i$
cannot be a descent of $I$ and if $c_i>c_{i+1}$ then $i$ has to be a descent
of $I$. This is precisely the definition of a cycloribbon.
All those simple modules are non-isomorphic since the eigenvalues
characterize them completely.
\end{proof}

We will prove in Section~\ref{induc0h} that all simple modules are of
dimension one, thus isomorphic to some $S_{[I,\c]}$ by following the
same pattern as in~\cite{BHT}, proving that all composition factors of the
modules induced from simple modules of the $0$-Hecke algebra are
one-dimensional.

Notice that if $[I,\c]$ is a cycloribbon, then $\eta_{[I,\c]}$ is an
eigenvector for all $T_i$, so that it can be written (since $B'$ is a basis
of $\AKSz$, see Proposition~\ref{ensBp}) as
$\eta_{I'} \sum_{\c'}{\alpha_{\c'} L_{\c'}}$, where $I'$ is the shape of the
anticycloribbon $\bijcyc([I,\c])$ and $\alpha_{\c'}$ are unknown coefficients.
The result is actually simpler:

\begin{proposition}
\label{ech-aks}
Let $[I,\c]$ be a cycloribbon. Then,
\begin{equation}
\eta_{[I,\c]} = \eta_{I'} L_{\overline{\c}},
\end{equation}
where $I'$ is the shape of the anticycloribbon $\bijcyc([I,\c])$.
\end{proposition}

\begin{proof}
Let $n:=|I|$ and let $\rho'$ be a color word involving $n$ distinct colors
such that
\begin{equation}
\left\{
\begin{array}{lcr}
\rho'_i<\rho'_j &\text{\ if\ } & \overline{\c}_i<\overline{\c}_j,\\
\rho'_i>\rho'_j &\text{\ if\ } & \overline{\c}_i>\overline{\c}_j.\\
\end{array}
\right.
\end{equation}
Let us define another color word $\rho$ as the permutation such that
\begin{equation}
\label{rhop}
\rho_i>\rho_j \text{\ iff\ }
\left\{
\begin{array}{lcl}
\rho'_i>\rho'_j &\text{and}& \c_i=\c_j, \text{\ or}\\
\rho'_i<\rho'_j &\text{and}& \c_i\not=\c_j.
\end{array}
\right.
\end{equation}
Thanks to Relation~(\ref{TL}), one proves by induction that
\begin{equation}
Y_\sigma(\rho')L_{\overline{c}}
= L_{\overline{c}\cdot\sigma^{-1}}Y_\sigma(\rho).
\end{equation}
Thanks to Lemmas~\ref{ech-hecke} and~\ref{omega-eta}, one then gets
\begin{equation}
Y_{\omega_n}(\rho') =\eta_{\overline{\Des(\rho')}}.
\end{equation}

Now, it is easy to build a permutation $\rho'$ satisfying these properties,
with descent composition equal to $\overline{I'}$ iff $I'$ is a composition
such that $[I',\c]$ is an anticycloribbon (the descent composition of $\rho'$
is, by definition of $\rho'$, already included in this set of compositions).
It then comes that
\begin{equation}
\eta_{I'} L_{\overline{c}}
= L_{\overline{c}\cdot\omega_n} \eta_{\overline{\Des(\rho)}}
= L_{\c} \eta_{\overline{\Des(\rho)}}.
\end{equation}
By definition of $\rho$, $[\overline{\Des(\rho)},\c]$ is the cycloribbon
$\bijcyc([I',\c])$.
\end{proof}

This result was not unexpected: it is clear that the generator of $S_{[I,\c]}$
has an expression as $\eta_{I'}L_{\c'}$, and from Lemma~\ref{omega-eta}
and Equation~(\ref{TL}), one can see that $\c'$ has to be
$\overline{\c}=\c\omega_n$.

\subsection{Induction of the simple $0$-Hecke modules}
\label{induc0h}

To describe the induction process, we need a partial order on
colored ribbons. Let $I$ be a composition and $\c=(c_1,\ldots,c_n)$.  The
covering relation of the order $\leq_I$ amounts to sort in increasing
order any two adjacent elements in the rows of $I$ or to sort in decreasing
order any two adjacent elements in the columns of $I$. For example, the
elements smaller than or equal to
\begin{equation}
\PetitTableau
T := \Tableau{2&1\\\ &1&3&3\\\ &\ &\ &4\\\ &\ &\ &3\\}
\end{equation}
are
\begin{equation}
\PetitTableau
\Tableau{2&1\\\ &1&3&3\\\ &\ &\ &4\\\ &\ &\ &3\\}
\qquad \qquad
\Tableau{1&2\\\ &1&3&3\\\ &\ &\ &4\\\ &\ &\ &3\\}
\qquad \qquad
\Tableau{2&1\\\ &1&3&4\\\ &\ &\ &3\\\ &\ &\ &3\\}
\qquad \qquad 
\Tableau{1&2\\\ &1&3&4\\\ &\ &\ &3\\\ &\ &\ &3\\}
\end{equation}

If $I$ is a composition of $n$, let $S_I:=H_n(0)\eta_I$ be the corresponding
simple module of $H_n(0)$ realized as a minimal left ideal in its left regular
representation (see~\cite{NCSF4}) and
\begin{equation}
M_I := S_I \uparrow_{H_n(0)}^{\AKSz}.
\end{equation}

Clearly, $M_I$ has dimension $r^n$ and admits $L_{\c} \eta_I$ as linear
basis, when $\c$ runs over color words.
For $\c\in C^n$, let $M_{\c,I}$ be the $\AKSz$-submodule of $M_I$
generated by $L_{\c}\eta_I$.

\begin{lemma}
\begin{equation}
M_{\c,I} \subseteq M_{\c',I} \quad\Longleftrightarrow\quad
\c\leq_I \c'.
\end{equation}
\end{lemma}

\begin{proof}
Let $i\in\{1,\ldots,n-1\}$.

\begin{itemize}
\item If $c_i<c_{i+1}$, and $T_i$ acts by $0$ on $\eta_I$, we get
$(1+T_i)L_{\c}\eta_I = L_{\c\sigma_i}T_i\eta_I=0$.
\item If $c_i<c_{i+1}$, and $T_i$ acts by $-1$ on $\eta_I$, we get
$(1+T_i)L_{\c}\eta_I = L_{\c\sigma_i}T_i\eta_I
=- L_{\c\sigma_i}\eta_I$, and so $-(1+T_i)$ sorts in decreasing order
$c_i$ and $c_{i+1}$.
\item If $c_i=c_{i+1}$, then $T_iL_{\c}\eta_I = L_{\c}T_i\eta_I$ so
the result is either $0$ or $-1$ times $L_{\c}\eta_I$.
\item If $c_i>c_{i+1}$, and $T_i$ acts by $-1$ on $\eta_I$, we get
$T_iL_{\c}\eta_I = L_{\c\sigma_i}(1+T_i)\eta_I = 0$.
\item If $c_i>c_{i+1}$, and $T_i$ acts by $0$ on $\eta_I$, we get
$T_iL_{\c}\eta_I = L_{\c\sigma_i}(1+T_i)\eta_I
= L_{\c\sigma_i}\eta_I$, and so $T_i$ sorts in increasing order
$c_i$ and $c_{i+1}$.
\end{itemize}
\end{proof}

We can now complete Theorem~\ref{simples} by proving that we have a complete
set of simple $\AKSz$-modules:

\begin{theorem}
The $S_{[I,\c]}$ form a complete family of simple modules of $\AKSz$.
\end{theorem}

\begin{proof}
It is sufficient to prove that any simple $\AKSz$-module $M$ appears as a
composition factor of some $M_I$ since by the previous lemma, $M_I$ admits a
composition series involving only one-dimensional modules.

Now, if $M$ is any simple $\AKSz$-module, let $S_I$ be a simple
$H_n(0)$-module occuring in the socle of $M\downarrow H_n(0)$. This means that
we have a monomorphism $i:S_i\hookrightarrow M$ of $H_n(0)$-modules. Inducing
to $\AKSz$, we get a non-zero $\AKSz$-morphism $f:M_I\to M$, which has to be
surjective since $M$ is simple.
So $M$ is a composition factor of some $M_I$ and all simple $\AKSz$-modules
are one-dimensional.
\end{proof}

\subsection{First Grothendieck ring}
\label{Groth-sec}

The induction product of two simple modules of $\AKSz$ can be worked out in a
way very similar to the case of $H_n(0)$-modules such as recalled in
Section~\ref{rthn}.
As in Proposition~\ref{ech-aks}, we represent these modules by their
anticycloribbons.
In the context of $\AKSz$, one has to extend to colored permutations what was
done with permutations in $H_n(0)$.

Let us consider two anticycloribbons $[I,\c]$ and $[I',\c']$.
The structure of the induction product $S_{[I,\c]}\widehat\otimes
S_{[I',\c']}$ of the corresponding simple modules is completely encoded in the
graph of the convolution of the associated colored permutations: as in
$H_n(0)$, there exists a basis $(b)$ of $S_{[I,\c]}\widehat\otimes
S_{[I',\c']}$ naturally labelled by the elements occuring in
the convolution, such that, either $T_i$, $-T_i$, or $1+T_i$ sends a given
basis element $b$ to another one, and $L_{\c} b$ is $b$ or $0$.
The unique one-dimensional quotient of the module generated by an element
$b$ is the one labelled by the anticycloribbon associated with its colored
permutation.
Let us describe more precisely this graph.

\begin{algorithm}
\label{induitalgo}
~

\noindent
\emph{Input}:  Two anticycloribbons $[I_1,\c_1]$ and $[I_2,\c_2]$.

\noindent
\emph{Output}:  A graph $G$.

Let $(\sigma_1,\c_1)$ and $(\sigma_2,\c_2)$ be the inverse of the
corresponding maximal colored permutations and let $\c=\c_1\concat\c_2$.

Let $G$ be the graph with vertices $V_w$ labelled by the elements of the
convolution $\sigma_1\convol\sigma_2$.
The edges of $G$ are labelled by the generators of our presentation of $\AKSz$
(see Proposition~\ref{notrepres}), that is, all Lagrange polynomials 
$L_{\bf d}$ and, for all $1\leq i\leq n-1$, the Hecke generator $T_i$, $-T_i$
or $1+T_i$.

The edge labelled by $L_{\bf d}$ on the vertex $V_\tau$ is a loop on
$V_\tau$ iff ${\bf d}=\c\cdot\tau^{-1}$.
Otherwise, $V_\tau$ is annihilated and the edge is not represented.
The edges labelled by the Hecke generators depend on the following cases:
\begin{itemize}
\item if the letters $\tau_i$ and $\tau_{i+1}$ come from the same factor of
the convolution, there is a loop labelled by $-T_i$ (resp. $1+T_i$) if $i$ is
(resp.  not) an anti-descent of $(\tau^{-1},\c\cdot\tau^{-1})$.
\item if $\tau_i$ comes from $\sigma_1$ and $\tau_{i+1}$ comes from
$\sigma_2$, there is a vertex from $V_\tau$ to $V_{\s_i\tau}$
labelled by $1+T_i$ (resp. $T_i$) if $i$ is (resp. not) an anti-descent of
$(\tau,\c)^{-1}$.
\item if $\tau_i$ comes from $\sigma_2$ and $\tau_{i+1}$ comes from
$\sigma_1$, there is a loop labelled by $-T_i$ (resp. $1+T_i$) if $i$ is
(resp. not) an anti-descent of $(\tau,\c)^{-1}$.
\end{itemize}
\end{algorithm}

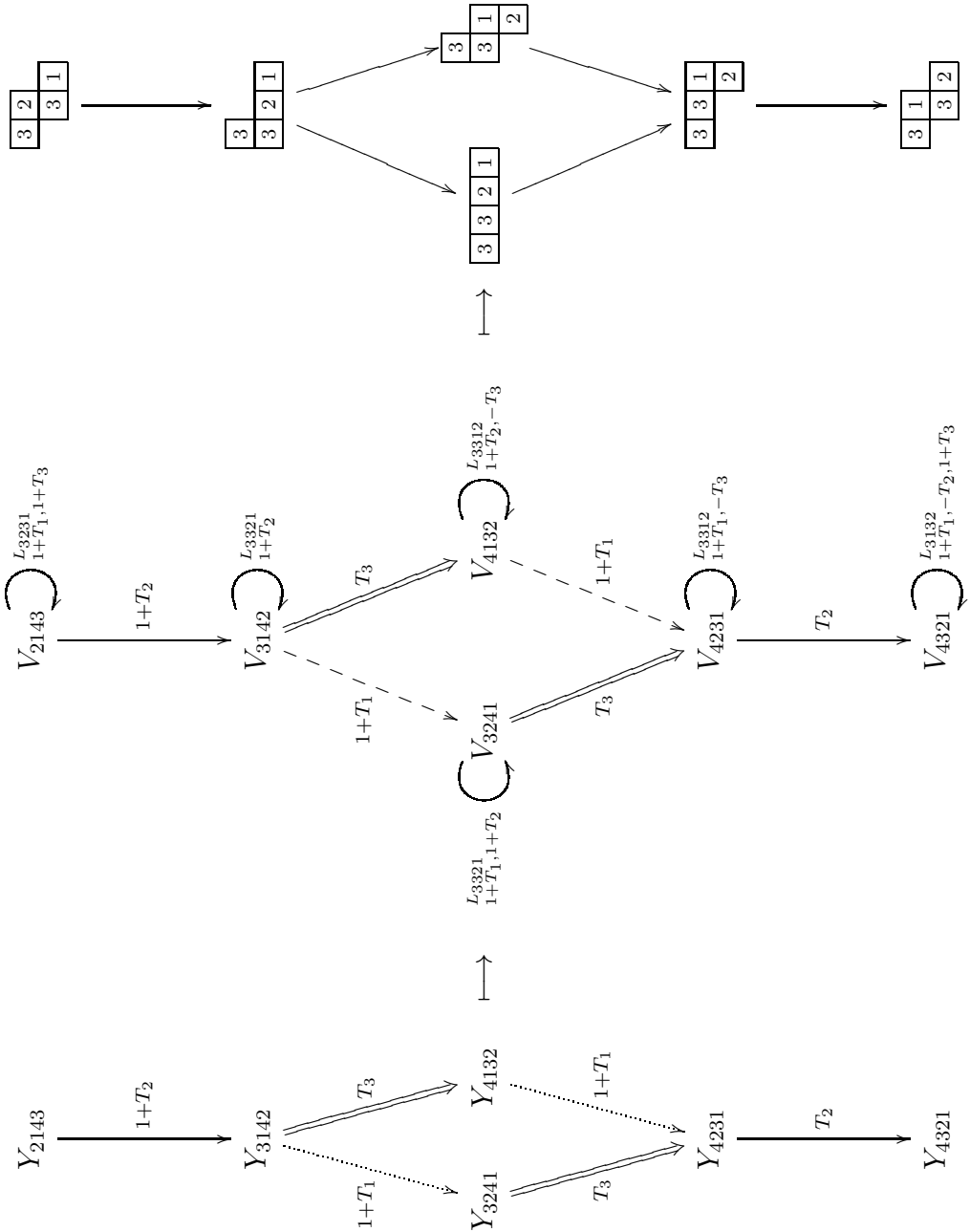
\begin{figure}[ht]
\PetitTableau
\ifdraft

\rotateleft{
\newdimen\vcadre\vcadre=0.2cm 
\newdimen\hcadre\hcadre=-0.2cm 
\def\GrTeXBox#1{\vbox{\vskip\vcadre\hbox{\hskip\hcadre%
      $\Tableau{#1\\}$%
   \hskip\hcadre}\vskip\vcadre}}
\def\GrTeXBoxb#1{\vbox{\vskip\vcadre\hbox{\hskip\hcadre%
   $#1$\hskip\hcadre}\vskip\vcadre}}

  \def\loopRight{\ar @`{[]+<1.3cm,1.2cm>,[]+<1.3cm,-1.2cm>}}
  \def\loopLeft{\ar @`{[]+<-1.3cm,1.2cm>,[]+<-1.3cm,-1.2cm>}}

$
\vcenter{\xymatrix@R=2.4cm@C=+0.4cm{
   & *{\GrTeXBoxb{Y_{2143}}} \ar @{->}[d]^{1+T_2} & \\
   & *{\GrTeXBoxb{Y_{3142}}}
   \ar @{..>}[dl]_{1+T_1} \ar @{=>}[dr]^{T_3}& \\ 
 *{\GrTeXBoxb{Y_{3241}}}\ar @{=>}[dr]_{T_3}
 & & *{\GrTeXBoxb{Y_{4132}}} \ar @{..>}[dl]^{1+T_1}\\
   & *{\GrTeXBoxb{Y_{4231}}} \ar @{->}[d]^{T_2}& \\
   & *{\GrTeXBoxb{Y_{4321}}} & \\
}}
\hcadre=0.0cm
\quad\longmapsto\quad
\vcenter{\xymatrix@R=2.4cm@C=+0.4cm{
    & *{\GrTeXBoxb{V_{2143}}} \loopRight^{{}^{L_{3231}}_{1+T_1, 1+T_3}} \ar @{->}[d]^{1+T_2} \\
    & *{\GrTeXBoxb{V_{3142}}} \loopRight^{{}^{L_{3321}}_{1+T_2}} \ar @{-->}[dl]_{1+T_1} \ar @{=>}[dr]^{T_3}\\
 *{\GrTeXBoxb{V_{3241}}} \loopLeft_{{}^{L_{3321}}_{1+T_1, 1+T_2}} \ar @{=>}[dr]_{T_3}
& & *{\GrTeXBoxb{V_{4132}}} \loopRight^{{}^{L_{3312}}_{1+T_2, -T_3}} \ar @{-->}[dl]^{1+T_1}\\
    & *{\GrTeXBoxb{V_{4231}}} \loopRight^{{}^{L_{3312}}_{1+T_1, -T_3}} \ar @{->}[d]^{T_2}\\
    & *{\GrTeXBoxb{V_{4321}}} \loopRight^{{}^{L_{3132}}_{1+T_1, -T_2, 1+T_3}} \\
}}
\hcadre=-0.2cm
\quad\longmapsto\quad
\vcenter{\xymatrix@R=1.8cm@C=+0.4cm{
   & *{\GrTeXBox{3 & 2\\\ & 3 & 1}} \ar @{->}[d] & \\
   & *{\GrTeXBox{3\\ 3 & 2 & 1}}
   \ar @{->}[dl] \ar @{->}[dr] & \\
     *{\GrTeXBox{3 & 3 & 2 & 1}}\ar @{->}[dr]
 & & *{\GrTeXBox{3\\ 3 & 1\\\ & 2}} \ar @{->}[dl] \\
   & *{\GrTeXBox{3& 3& 1\\\ &\ &2}} \ar @{->}[d] & \\
   & *{\GrTeXBox{3& 1\\\ & 3 & 2}} & \\
}}
$
}
\fi
\caption{\label{ex-induc}
Induction product of the two simple modules
$S_{[(1,1),(23)]}$ and $S_{[(1,1),(13)]}$ of $\AKSzgen{2}{2}$ to
$\AKSzgen{4}{2}$.
The first graph is the restriction of the left weak order to the convolution
$21\convol21$.
In the second graph, we show the graph built by Algorithm~\ref{induitalgo},
and the third graph represents the composition factors corresponding to each
basis element.}
\end{figure}

For example, the second graph of Figure~\ref{ex-induc} is the graph
associated with the anticycloribbons $[(1,1),23]$ and $[(1,1),13]$.

\begin{theorem}
\label{th-induit}
Let $[I_1,\c_1]$ and $[I_2,\c_2]$ be two anticycloribbons.
Let $\CC V$ the vector space with basis $V$ indexed by the vertices of the
graph of Algorithm~\ref{induitalgo}. Then the relations $\rho V_i=V_j$
whenever $V_i\xrightarrow{\rho}V_j$ is a presentation of the induced
module $S_{[I_1,\c_1]}\widehat{\otimes} S_{[I_2,\c_2]}$.
Its composition factors are given by the shapes of the anticycloribbons
corresponding to the colored permutations $(\tau,\c)^{-1}$,
where $\tau$ runs over the labellings of the vertices $V$.
\end{theorem}

For example, Figure~\ref{ex-induc} presents the induction product of the two
simple modules $S_{[(1,1),23]}$ and $S_{[(1,1),13]}$ (written as
anticycloribbons) of $\AKSzgen{2}{2}$ to $\AKSzgen{4}{2}$.
In particular, we get the following composition factors (written as
anticycloribbons):
\begin{equation}
\PetitTableau
\Tableau{3&2\\\ &3&1\\}
\qquad
\Tableau{3\\ 3&2&1\\}
\qquad
\Tableau{3&3&2&1\\}
\qquad
\Tableau{3\\3&1\\\ &2\\}
\qquad
\Tableau{3&3&1\\\ &\ &2\\}
\qquad
\Tableau{3&1\\\ &3&2\\}
\end{equation}
Note that this construction gives an effective algorithm to compute the
composition factors of the induction product of two simple modules. One
recovers the same shapes as in Figure~\ref{0-Hecke-graphe}, that is, the
induction of the simple modules $(2)$ and $(2)$ of $H_2(0)$ in $H_4(0)$.
Moreover, one recovers these shapes in the same order as in
Figure~\ref{fig-yb}. This property comes from Proposition~\ref{ech-aks}, since
any basis element $V_{(\sigma,\c)}$ can be written either in the
form $L_{\d}T$ or $T'L_{\d'}$, where $T$ and $T'$ are in $H_n(0)$.

\begin{proof}
%
Let us consider two simple modules $M=S_{[I_1,\c_1]}$, $N=S_{[I_2,\c_2]}$ of
$\AKSzgen{m}{r}$ and $\AKSzgen{n}{r}$ and let $\sigma_1$ and $\sigma_2$ be
the maximal colored permutations associated with their anticycloribbons.
Let $\c=\c_1\concat\c_2$ be the concatenation of the color words.

The induced module $M\widehat{\otimes} N$ of $\AKSzgen{m+n}{r}$ is by
definition
\begin{equation}
\AKSzgen{m+n}{r} \bigotimes_{\AKSzgen{m}{r}\otimes\AKSzgen{n}{r}}
\left(M\otimes_{\CC} N\right).
\end{equation}
Since every $L_{\c}$ belongs to $\AKSzgen{m}{r}\otimes\AKSzgen{n}{r}$, the
induced module is isomorphic to
$H_{m+n}(0)\otimes_{H_m(0)\otimes H_n(0)} (M\otimes N)$ as a vector space and
even as a $H_{m+n}(0)$-module. In particular, its dimension is
$\binom{m+n}{n}$.

Now, as a $H_{m+n}(0)$-module, thanks to Lemma~\ref{yb-lemme},
$M\widehat{\otimes} N$ is described by the restriction $G_1$ of the graph of
$\SG_{m+n}$ to $\sigma_1\convol\sigma_2$, with the
Yang-Baxter elements determined by the color word $\c$.
Let $G$ be the graph associated with $M$ and $N$ by
Algorithm~\ref{induitalgo}.
The edges of $G$ and $G_1$ are the same since both are edges in
$\sigma_1\convol\sigma_2$ with the same constraints.
The loops of $G$ and $G_1$ labelled by the Hecke generators are the same
thanks to Algorithm~\ref{induitalgo} and Corollary~\ref{yb-bases}.
So $G$ satisfies the relations of the $0$-Hecke algebra,
relations~(\ref{ti1}), (\ref{ti2}), and~(\ref{ti3}).
It trivially satisfies Relations~(\ref{LL}) and it satisfies
Relations~(\ref{TL}), since each vertex of the graph can be written as
some Lagrange polynomial $L_{\c}$ multiplied by an element of $H_n(0)$.
So it satisfies all the relations of our presentation of the
Ariki-Koike-Shoji algebra at $q=0$ (Proposition~\ref{notrepres}).
Since it has the correct dimension and the correct generator, this graph
encodes the structure of the induced module $M\widehat{\otimes} N$.

Finally, the composition factors are the anticycloribbons associated with the
colored permutations $(\tau,\c)^{-1}$: the only non-zero $L$ is
the right one and the action of $T_i$ is the one encoded by the shape of the
anticycloribbon, thanks to Note~\ref{simpleYB}.
\end{proof}

\subsection{Restriction of simple modules}
\label{restr-sec}

\begin{proposition}
Let $[I,\c]$ be a cycloribbon, regarded as a filling of the ribbon diagram of
$I$. Let $I'$ and $I''$ be the compositions whose diagrams are formed
respectively of the first $k$ cells and the last $n-k$ cells of $I$ and
$\c'=(c_1,\ldots,c_k)$, $\c''=(c_{k+1},\ldots,c_n)$. Then
$[I',\c']$ and $[I'',\c'']$ are again cycloribbons and the restriction of
$S_{[I,\c]}$ to $\AKSzgen{k}{r}\otimes\AKSzgen{n-k}{r}$ is
$S_{[I',\c']}\otimes S_{[I'',\c'']}$.
\end{proposition}

\begin{figure}[ht]
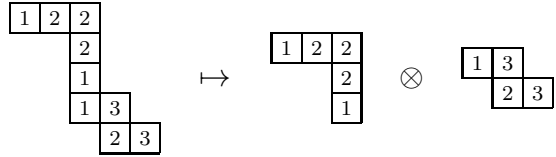

$$
\PetitTableau
\Tableau{1&2&2\\\ &\ &2\\\ &\ &1\\\ &\ &1&3 \\\ &\ &\ &2&3\\}
\quad \mapsto \quad
\Tableau{1&2&2\\\ &\ &2\\\ &\ &1\\}
\quad \otimes \quad
\Tableau{1&3 \\\ &2&3\\}
$$
\caption{\label{restr-simp}The restriction of the simple module represented on
the left by a cycloribbon to $\AKSzgen{5}{r}\otimes\AKSzgen{4}{r}$ gives the
simple module on the right also represented by cycloribbons.}
\end{figure}

For example, Figure~\ref{restr-simp} shows that the restriction of the first
module
$S_{[31122,122211323]}$ to $\AKSzgen{5}{r}\otimes\AKSzgen{4}{r}$ is
$S_{[311,12221]}\otimes S_{[22,1323]}$.

\begin{proof}
This follows from the obvious graphical description of the restriction of
simple $H_n(0)$-modules, which consists in cutting the ribbon diagram,
as was done before on cycloribbons, and from the commutation
relations between $L$ and $T$.
\end{proof}

\subsection{The quasi-symmetric Mantaci-Reutenauer algebra}

A \emph{colored composition} is a pair
$(I,u)=((i_1,\ldots,i_m),(u_1,\ldots,u_m))$ formed of a composition $I$ and a
color word of the same length. 

\begin{figure}[ht]
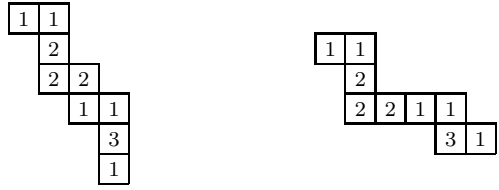

\PetitTableau
$$
\Tableau{1&1 \\\ &2\\\ &2&2 \\\ &\ &1&1 \\\ &\ &\ & 3 \\\ &\ &\ & 1\\}
\qquad
\qquad
\qquad
\Tableau{1&1 \\\ &2\\\ &2&2&1&1 \\\ &\ &\ &\ & 3 & 1\\}
$$
\caption{\label{cc-acr}A colored composition and the corresponding
anticycloribbon.}
\end{figure}

There is a bijection between colored compositions and anticycloribbons:
starting with a colored composition, one rebuilds an anticycloribbon by
separating adjacent blocks of the colored composition with different colors
and gluing them back in the only possible way according to the criterion of
being an anticycloribbon. Conversely, one separates two adjacent blocks of
different colors and glue them one below the other.
For example, Figure~\ref{cc-acr} shows a colored composition represented as a
filling of a ribbon ($u_i$ is written in row $i$) and its corresponding
anticycloribbon.

Let $X^{(i)}=\{x^{(i)}_k | k\geq1\}$ for $1\leq i\leq r$ be $r$ infinite
linearly ordered sets of commuting indeterminates. The \emph{monochromatic
monomial quasi-symmetric functions} labelled by a colored composition
$(I,u)$ is
\begin{equation}
M_{(I,u)}(X) = \sum_{j_1<\cdots<j_m}
(x_{j_1}^{(u_1)})^{i_1} \cdots (x_{j_m}^{(u_m)})^{i_m}.
\end{equation}

It is known (\cite{NT-coul} Proposition 5.2), that the $M_{I,u}$ span a
subalgebra $\QMR^{(r)}$ of $\CC[X^{(1)},\ldots,X^{(r)}]$.
This algebra is also a subalgebra of the algebra of level $r$ quasi-symmetric
functions defined in~\cite{Poi} (see also~\cite{NT-coul}).

We will prove that the Grothendieck ring of the tower of $\AKSz$ algebras is
isomorphic to $\QMR^{(r)}$.

Define an order on colored compositions as follows: let
\begin{equation*}
(I,u)= ((i_1,\ldots,i_m),(u_1,\ldots,u_m)) \quad\text{\ and\ }\quad
(J,v) = ((j_1,\ldots,j_p),(v_1,\ldots,v_p))
\end{equation*}
two colored compositions. Then $(I,u)$ is \emph{finer} than $(J,v)$ if
there exists a sequence $(l_0=0,l_1,\ldots,l_p=m)$ such that
for any integer $k$,
\begin{equation}
j_k=i_{l_{k-1}+1}+\cdots+i_{l_k} \quad\text{\ and\ }\quad
v_k=u_{l_{k-1}+1}=\cdots=u_{l_k}.
\end{equation}
For example, the compositions fatter than $((1,1,3,2),(2,1,1,2))$ are
\begin{equation}
((1,1,3,2),(2,1,1,2)) \quad \text{\ and\ }\quad ((1,4,2),(2,1,2)).
\end{equation}

This allows us to define the \emph{monochromatic quasi-ribbon functions}
$F_{(I,u)}$ by
\begin{equation}
F_{(I,u)} = \sum_{{(I',u')}\finer {(I,u)}} M_{(I',u')}.
\end{equation}
Notice that this last description of the order $\finer$ is reminiscent of the
order $\finer'$ on descent sets used in the context of quasi-symmetric
functions and non-commutative symmetric functions: more precisely, one gets
the usual order when considering compositions colored with
only one color.

\begin{lemma}[\cite{BH,NT-coul}]
The $F_{(I,u)}$, where $(I,u)$ runs over colored compositions,
span $\QMR^{(r)}$. This algebra is isomorphic to the dual of the
Mantaci-Reutenauer algebra $\MR^{(r)}$ defined in~\cite{MR}.
\qed
\end{lemma}

Let 
\begin{equation}
{\mathcal G}^{(r)} := \bigoplus_{n\geq 0} G_0(\AKSz)
\end{equation}
be the Grothendieck ring of the tower of $\AKSz$ algebras. Define a
characteristic map
\begin{equation}
\begin{array}{rcl}
ch: & {\mathcal G}^{(r)} &\to \QMR^{(r)}\\
& [S_{[I,\c]}] &\mapsto F_{[I,\c]}.
\end{array}
\end{equation}
Comparing the descriptions of $S_{[I',\c']}\widehat\otimes S_{[I'',\c'']}$ and
of $S_{[I,\c]}\downarrow$ obtained in Sections~\ref{Groth-sec}
and~\ref{restr-sec} with the formulas for product and coproduct of the $F$
basis in $\QMR^{(r)}$, we obtain:

\begin{theorem}
$ch$ is an isomorphism of Hopf algebras.
\qed
\end{theorem}


\newpage
\section{Projective modules}

\subsection{The Grothendieck ring}

The previous results imply, by duality, a description of the
Grothendieck ring of the category of projective $\AKSz$-modules.
Indeed, we already know that $\AKSz$ is self injective. Moreover, it is easy
to see that $\AKSz$ is a free module and hence projective over any parabolic
subalgebra
\begin{equation}
\AKSzgen{n_1}{r} \otimes\cdots\otimes \AKSzgen{n_k}{r}.
\end{equation}
These conditions are sufficient to ensure that, as a Hopf algebra,

\begin{equation}
{\mathcal K} = \bigoplus_{n\geq0} K_0(\AKSz) \simeq \MR^{(r)}
\end{equation}
is the graded dual of ${\mathcal G}$, and thus is isomorphic to the
Mantaci-Reutenauer algebra. Under this isomorphism, the classes of
indecomposable projective modules are mapped to the subfamily of the dual
basis $F_{[I,\c]}$ of Poirier quasi-symmetric functions labelled by colored
descent sets, or cycloribbons.

Recall that the Hopf algebra $\MR^{(r)}$ can be defined (see~\cite{NT-coul})
as the free associative algebra over symbols
$(S_{j}^{(i)})_{j\geq1 ; 1\leq i\leq r}$, graded by $\deg S_j^{(i)}=j$, and
with coproduct
\begin{equation}
\Delta S_n^{(k)} = \sum_{i=0}^n S_i^{(k)} \otimes S_j^{(k)}.
\end{equation}
For a colored composition $(I,u)$ as above, one defines
\begin{equation}
S^{(I,u)} := S_{i_1}^{(u_1)} \cdots S_{i_p}^{(u_p)}.
\end{equation}
Clearly, the $S^{(I,u)}$ form a linear basis of $\MR^{(r)}$.

The \emph{monochromatic colored ribbon basis} $R_{(I,u)}$ of $\MR^{(r)}$ can now be defined
by the condition

\begin{equation}
S^{(I,u)} =: \sum_{(J,v)\leq (I,u)} R_{(J,v)}.
\end{equation}

Let $(I,u) = (i_1,\ldots, i_p ; u_1,\ldots, u_p)$ and
$(J,v)= (j_1,\ldots, j_q;v_1,\ldots, v_q)$ be two colored compositions.
We set
\begin{equation}
(I,u) \concat (J,v) := (I\concat J,u\concat v),
\end{equation}
where $a\concat b$ denotes concatenation.

Moreover, if $u_p=v_1$, we set
\begin{equation}
(I,u)\triangleright (J,v) :=
( (i_1,\ldots, i_{p-1},(i_{p}+j_1),j_2,\ldots,j_q );
 (u_1,\ldots, u_p,v_2,\ldots, v_q)).
\end{equation}

The colored ribbons satisfy the very simple multiplication rule:
\begin{equation}
\label{prodRcol}
R_{(I,u)} R_{(J,v)} = R_{(I,u)\concat (J,v)} +
\left\{
\begin{array}{ll}
R_{(I,v)\triangleright(J,v)} & \text{if}\ u_p=v_1,\\
0 & \text{if}\ u_p\not=v_1.
\end{array}
\right.
\end{equation}

Let $[K,\c]$ be an anticycloribbon. Let $P_{[K,\c]}$ be the
indecomposable projective module whose unique simple quotient is the simple
module labelled by $\bijcyc([K,\c])$ and let $(I,u)$ be the corresponding
colored composition.
Summarizing this discussion, we have

\begin{theorem}
\label{K2MR}
The map
\begin{equation}
\begin{array}{rcl}
{\bf ch}: & {\mathcal K}     &\to \MR^{(r)} \\
          &[P_{[K,\c]}] &\mapsto R_{(I,u)}
\end{array}
\end{equation}
is an isomorphism of Hopf algebras.
\qed
\end{theorem}

The main interest of the labelling by colored compositions is that it allows
immediate reading of some important information. For example, it follows from
Theorem~\ref{K2MR} that the products of complete functions $S^{(I,u)}$ are the
characteristics of the projective $\AKSz$-modules obtained as induction
products of the one-dimensional projective $\AKSzgen{m}{r}$-modules on which
all $T_i$ act by $0$ and all $\xi_j$ by the same eigenvalue $u_i$.

We see that, as in the case of $H_n(0)$, each indecomposable projective
$\AKSz$-module occurs as a direct summand of such an induced module, and that
the direct sum decomposition is given by the anti-refinement order.
For example, writing a bar over the parts of the composition corresponding to
color $2$ (and nothing over the parts corresponding to color $1$), the
identity
\begin{equation}
S^{(2\overline{1}\overline{2}13)} =
R_{ (2\overline{1}\overline{2}13)} +
R_{(2\overline{3}13)} +
R_{(2\overline{1}\overline{2}4)} +
R_{(2\overline{3}4)}
\end{equation}
indicates which indecomposable projective direct summands compose the
projective $\AKS_{9,2}(0)$-module defined as the outer tensor product
\begin{equation}
S_2 \,\widehat\otimes\, S_{\overline{1}} \,\widehat\otimes\, S_{\overline{2}}
\,\widehat\otimes\, S_1 \,\widehat\otimes\, S_3.
\end{equation}

Since by multiplying a Lagrange element $L_\c$ of $\AKSz$ by any element of
$\AKSz$, one can only obtain an element of the form $L_\d T$ where $\d$ has
same evaluation as $\c$ and $T$ is any Hecke element, we have a direct sum
decomposition of $\AKSz$ into two-sided ideals:
\begin{equation}
\AKSz = \bigoplus_{e} H_{n,r}^{e},
\end{equation}
where $e$ is any evaluation of $C^n$ and
$H_{n,r}^{e} := \AKSz \sum_{\c\in e} L_\c$.
Actually all such sums are central idempotents of $\AKSz$.
The ideals $H_{n,r}^e$ where $e=(i,i,\ldots,i)$ are said to be
\emph{monochromatic} of color $i$.
The monochromatic ideals are isomorphic to $H_n(0)$ as algebras, and their
indecomposable projective modules $P_{[I,\c]}$ are obtained by defining the
action of all $L_\c$ by zero except $L_{(i^n)}$ acting by one on the
$H_n(0)$-modules $P_I$.
Such modules are also called monochromatic. Note that one can realize
$P_{[I,i^k]}$ in the left regular representation as
\begin{equation}
P_{[I,i^k]} = \AKSz L_{i^n} \nu_I.
\end{equation}

It follows from Theorem~\ref{K2MR} and Formula~(\ref{prodRcol}) that
\begin{corollary}
\label{indmono}
Every indecomposable projective module $P_{[I,\c]}$ of $\AKSz$ is obtained as
an induction product of monochromatic modules
\begin{equation}
\label{proji}
P_{[I,\c]} = P_{[I_1,\c_1]}\widehat\otimes \cdots \widehat\otimes
P_{[I_k,\c_k]},
\end{equation}
where $I=(I_1,\ldots,I_k)$ and $I_j$ is a composition,
$\c=(\c_1,\ldots,\c_k)$ and $\c_j$ is a monochromatic word.
\qed
\end{corollary}

Corollary~\ref{indmono} and Theorem~\ref{induit-shuf} implies that
$P_{[I,\c]}$ is a semi-combinatorial module. Indeed, it is possible to choose
for the underlying $H_n(0)$-module
$P_{I_1}\widehat\otimes \cdots \widehat\otimes P_{I_k}$ a Yang-Baxter graph
such that all vertices are common eigenvectors of the Lagrange polynomials
by making the right choice on the color vector $\d$: if $c_i<c_j$ then
$d_i<d_j$ and if $c_i>c_j$ then $d_i>d_j$.

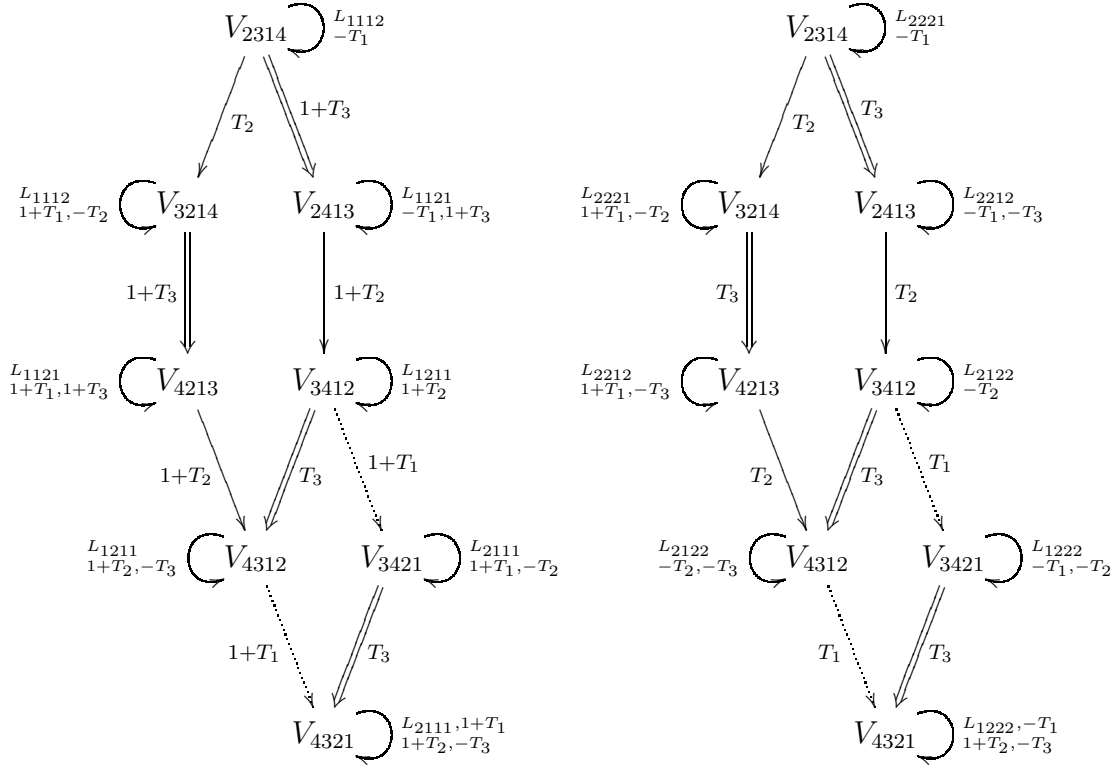
\begin{figure}[ht]
\PetitTableau
\ifdraft

{
\newdimen\vcadre\vcadre=0.2cm 
\newdimen\hcadre\hcadre=0.0cm 
\def\GrTeXBox#1{\vbox{\vskip\vcadre\hbox{\hskip\hcadre%
      $\Tableau{#1\\}$%
   \hskip\hcadre}\vskip\vcadre}}
\def\GrTeXBoxb#1{\vbox{\vskip\vcadre\hbox{\hskip\hcadre%
   $#1$\hskip\hcadre}\vskip\vcadre}}
  \def\loopRight{\ar @`{[]+<1.2cm,1.1cm>,[]+<1.2cm,-1.1cm>}}
  \def\loopLeft{\ar @`{[]+<-1.2cm,1.1cm>,[]+<-1.2cm,-1.1cm>}}
  \def\loopRightb{\ar @`{[]+<1.5cm,1.1cm>,[]+<1.5cm,-1.1cm>}}
  \def\loopLeftb{\ar @`{[]+<-1.8cm,1.1cm>,[]+<-1.8cm,-1.1cm>}}

$
\vcenter{\xymatrix@R=1.6cm@C=+0.05cm{
   & *{\GrTeXBoxb{V_{2314}}} \loopRight^{{}^{L_{1112}}_{-T_1}}
     \ar @{->}[dl]^{T_2} \ar @{=>}[dr]^{1+T_3}& \\
 *{\GrTeXBoxb{V_{3214}}}\loopLeft_{{}^{L_{1112}}_{1+T_1,-T_2}}
  \ar @{=>}[d]_{1+T_3}
 & & *{\GrTeXBoxb{V_{2413}}}\loopRight^{{}^{L_{1121}}_{-T_1,1+T_3}}
  \ar @{->}[d]^{1+T_2}\\
 *{\GrTeXBoxb{V_{4213}}}\loopLeft_{{}^{L_{1121}}_{1+T_1,1+T_3}}
  \ar @{->}[dr]_{1+T_2}
 & & *{\GrTeXBoxb{V_{3412}}}\loopRight^{{}^{L_{1211}}_{1+T_2}}
  \ar @{=>}[dl]^{T_3} \ar @{..>}[dr]^{1+T_1}\\
   & *{\GrTeXBoxb{V_{4312}}}\loopLeft_{{}^{L_{1211}}_{1+T_2,-T_3}}
  \ar @{..>}[dr]_{1+T_1}
   & & *{\GrTeXBoxb{V_{3421}}}\loopRight^{{}^{L_{2111}}_{1+T_1,-T_2}} 
  \ar @{=>}[dl]^{T_3} \\
 & & *{\GrTeXBoxb{V_{4321}}}\loopRight^{{}^{L_{2111},1+T_1}_{1+T_2,-T_3}} & \\
}}
\vcenter{\xymatrix@R=1.6cm@C=+0.05cm{
   & *{\GrTeXBoxb{V_{2314}}} \loopRight^{{}^{L_{2221}}_{-T_1}}
  \ar @{->}[dl]^{T_2} \ar @{=>}[dr]^{T_3}& \\
 *{\GrTeXBoxb{V_{3214}}}\loopLeft_{{}^{L_{2221}}_{1+T_1,-T_2}}
  \ar @{=>}[d]_{T_3}
 & & *{\GrTeXBoxb{V_{2413}}}\loopRight^{{}^{L_{2212}}_{-T_1,-T_3}}
  \ar @{->}[d]^{T_2}\\
 *{\GrTeXBoxb{V_{4213}}}\loopLeft_{{}^{L_{2212}}_{1+T_1,-T_3}}
  \ar @{->}[dr]_{T_2}
 & & *{\GrTeXBoxb{V_{3412}}}\loopRight^{{}^{L_{2122}}_{-T_2}}
  \ar @{=>}[dl]^{T_3} \ar @{..>}[dr]^{T_1}\\
   & *{\GrTeXBoxb{V_{4312}}}\loopLeft_{{}^{L_{2122}}_{-T_2,-T_3}}
  \ar @{..>}[dr]_{T_1}
   & & *{\GrTeXBoxb{V_{3421}}}\loopRight^{{}^{L_{1222}}_{-T_1,-T_2}}
  \ar @{=>}[dl]^{T_3} \\
 & & *{\GrTeXBoxb{V_{4321}}}\loopRight^{{}^{L_{1222},-T_1}_{1+T_2,-T_3}} & \\
}}
$
}
\fi
\caption{\label{ind-proj-aks}
Induction product of the two indecomposable projective modules
$P_{[(2,1),(11)]}\widehat\otimes P_{[(1),(2)]}$ (left module).
Induction product of the two indecomposable projective modules
$P_{[(2,1),(22)]}\widehat\otimes P_{[(1),(1)]}$ (right module).
}
\end{figure}

Figure~\ref{ind-proj-aks} shows the two examples of induction product
$P_{[(2,1),(1)]}\widehat\otimes P_{[(1),(2)]}$ and
$P_{[(2,1),(2)]}\widehat\otimes P_{[(1),(1)]}$. On the first induced module,
the color vector has to be in the convolution of $312$ and $1$ and have its
last element greater than all the others, and so is $3124$.
On the second induced module, the color vector has to be in the convolution of
$312$ and $1$ and have its last element smaller than all the others, and so is
$4231$. One can compare this figure with Figure~\ref{fig-projind}.

Their restrictions to $H_n(0)$ (and hence, their dimensions) can be computed
by means of the following result:

\begin{corollary}
\label{mr-ncsf}
The homomorphism of Hopf algebras
\begin{equation}
\label{homom}
\begin{array}{rcl}
\pi: &\MR^{(r)}    & \to \NCSF \\
     & S_{j}^{(i)} & \mapsto S_j
\end{array}
\end{equation}
maps the class of a projective $\AKSz$-module to the class of its restriction
to $H_n(0)$.
\end{corollary}

\begin{proof}
The previous considerations show in particular that the restriction of a
mono\-chromatic $\AKSz$-module $P_{[I,i^n]}$ to $H_n(0)$ is isomorphic to the
indecomposable projective $H_n(0)$-module $P_I$.
Hence, the restriction of any indecomposable projective module $P_I$ given as
in Formula~(\ref{proji}) is
\begin{equation}
P_{I_1}\widehat\otimes \cdots \widehat\otimes P_{I_k},
\end{equation}
so that the restriction map $\pi$ is given by
\begin{equation}
\pi(R_{[I,\c]}) = R_{I_1} \cdots R_{I_k},
\end{equation}
which is equivalent to Formula~(\ref{homom}).
\end{proof}

Continuing the previous example, we see that the restriction of
$P_{2\overline{1}\overline{2}13}$ to $H_9(0)$ is given by
\begin{equation}
\pi(R_{2\overline{1}\overline{2}13}) =
\pi(R_2R_{\overline{12}}R_{13}) =
R_2R_{12} R_{13} =
R_{21213} + R_{2133} + R_{3213} + R_{333}.
\end{equation}

Dually, one can describe the induction of projective $H_n(0)$-modules
to $\AKSz$. The next result is a consequence of Corollary~\ref{mr-ncsf}.

\begin{corollary}
Let $I$ be a composition of $n$ and let $N_I$ be the $\AKSz$-module induced by
the indecomposable projective $H_n(0)$-module $P_I$.
Then
\begin{equation}
N_I \simeq \bigoplus P_{[I,\c]},
\end{equation}
where the sum runs over all the anticycloribbons of shape $I$.
\qed
\end{corollary}

For example, let us complete the case $I=(2,1)$ with two colors. The
following five anticycloribbons appear in the induction of $P_I$, with
respective dimensions $3$, $6$, $3$, $2$, and $2$:
\begin{equation}
\PetitTableau
\Tableau{2&1\\\ &1\\}
\qquad
\Tableau{2&1\\\ &2\\}
\qquad
\Tableau{1&1\\\ &2\\}
\qquad
\Tableau{2&2\\\ &2\\}
\qquad
\Tableau{1&1\\\ &1\\}
\end{equation}

\subsection{Cartan invariants and decomposition numbers}

Finally, we can describe the Cartan invariants and the decomposition matrices
by means of the maps
$$
\xymatrix{
Sym^{(r)} \ar[rr]^d &                             & QMR^{(r)} \\
                & {\MR^{(r)}}\ar[lu]_e \ar[ru]_c & \\
}
$$

\begin{equation}
\begin{array}{rcl}
e: &\MR^{(r)} &\to Sym^{(r)} = (Sym)^{\otimes r} \simeq Sym(X_1,\ldots,X_r),\\
   & S_j^{(i)} &\mapsto h_j(X_i)
\end{array}
\end{equation}
and
\begin{equation}
\begin{array}{rcl}
d:& Sym^{(r)} &\hookrightarrow \QMR^{(r)}\\
& h_j(X_i) &\mapsto F_{[(j),i^j]}.
\end{array}
\end{equation}

Then the Cartan map is $c=d\circ e$, and the entry
$c_{[I,\c], [J,{\bf d}]}$ of the Cartan matrix giving the multiplicity
of the simple module $S_{[J,{\bf d}]}$ as a composition factor of the
indecomposable projective module $P_{[I,\c]}$ is equal to the
coefficient of $F_{[J,{\bf d}]}$ in $c(R_{[I,\c]})$.


\medskip
For an $r$-partition $\rlambda=(\lambda^{(1)},\ldots,\lambda^{(r)})$,
let $Z_\rlambda(q)$ be the irreducible module of the generic algebra
$\AKSq$ as constructed in~\cite{SS}, Theorem 4.1.
With our normalization, this module is defined for arbitrary values of $q$,
including $q=0$. For an $r$-colored composition $(I,u)$ of $n$, let
\begin{equation}
\label{eq100}
d_{\rlambda,(I,u)} := [ S_{(I,u)} : Z_\rlambda(0)]
\end{equation}
be the decomposition numbers.

\begin{theorem}
Let $r_{(I,u)}=e(R_{(I,u)})$ be the commutative image of $R_{(I,u)}$ in
$Sym^{(r)}$.
Let $S_\rlambda = s_{\lambda^{(1)}}(X_1)\cdots s_{\lambda^{(r)}}(X_r)$
as in~\cite{Sho}. Then the multiplicity of the simple module $S_{(I,u)}$ as a
composition factor of $Z_\rlambda(0)$ is equal to the coefficient of
$S_\rlambda$ in $r_{(I,u)}$, that is
\begin{equation}
d_{\rlambda,(I,u)} := \langle S_\rlambda, r_{(I,u)} \rangle
\end{equation}
where $\langle,\rangle$ is the scalar product on $Sym^{(r)}$ for which the
basis $(S_\rlambda)$ is orthonormal.
\end{theorem}

\begin{proof}
This follows from Shoji's character formula (Theorem 6.14 of~\cite{Sho}).
We first need to reformulate it on another set of elements. Instead of the
$a_{\mu}$ of Shoji, we shall make use of elements $b_{(J,v)}$
labelled by colored compositions, and defined by
\begin{equation}
b_{(J,v)} := L_{\c(J,v)} T_{w(J)}
\end{equation}
where $\c(J,v) = (v_1^{j_1}\cdots v_s^{j_s})$ and
$w(J) = \gamma_{j_1}\times\cdots\times\gamma_{j_s}\in
\SG_{j_1}\times\cdots\times\SG_{j_s}$, each $\gamma_k=s_{k-1}\cdots s_1$ being
a $k$-cycle.
With each $b_{(J,v)}$, we associate a ``non-commutative cycle index''
\begin{equation}
{\bf Z}(b_{(J,v)}) := (q-1)^{-l(J)}S^{(J,v)}((q-1)A)\in \MR_n^{(r)}.
\end{equation}
Then, Shoji's formula (6.14.1) is equivalent to
\begin{equation}
\chi_q^\rlambda(b_{(J,v)}) = \langle S_\rlambda, {\bf Z}(b_{(J,v)})\rangle
\end{equation}
where $S_\rlambda$ is interpreted as an element of $QMR^{(r)}$.

Specializing this at $q=0$ gives
\begin{equation}
\chi_0^\rlambda(b_{(J,v)}) =
\langle S_\rlambda, (-1)^{n-l(I)}\Lambda^{(J,v)}(A) \rangle,
\end{equation}
where $\Lambda^{(J,v)}(A) = \Lambda_{j_1}^{(v_1)}\cdots\Lambda_{j_s}^{(v_s)}$,
and the colored elementary functions $\Lambda_j^{(i)}$ are
\begin{equation}
\sum_{j}(-1)^j \Lambda_j^{(i)} = \left(\sum_k S_k^{(i)}\right)^{-1}.
\end{equation}
On another hand, using the two expressions
\begin{equation}
\eta_{(I,u)} = L_{\c'} \eta_{I'} = \eta_{I''} L_{\c''}
\end{equation}
of the generator of $S_{(I,u)}$, one checks that the eigenvalue of
$b_{(J,v)}$ on $\eta_{(I,u)}$ is
\begin{equation}
\langle F_{(I,u)} , {\bf Z}(b_{(J,v)}) \rangle.
\end{equation}
Therefore, $\chi_0^\rlambda(b_{(J,v)})$ is also given by
\begin{equation}
\chi_0^\rlambda(b_{(J,v)}) =
\sum_{(I,u)} d_{\rlambda,(I,u)} \langle F_{(I,u)} , {\bf Z}(b_{(J,v)})
\rangle.
\end{equation}
The ${\bf Z}(b_{(J,v)})$ being linearly independent, one can conclude that
\begin{equation}
\sum_{(I,u)} d_{\rlambda,(I,u)} F_{(I,u)} = S_\rlambda,
\end{equation}
which is equivalent to Equation~(\ref{eq100}) since $r_{(I,u)}$ is the image
by $e$ of $R_{(I,u)}$ which is itself the dual basis of $F_{(I,u)}$.
\end{proof}

\newpage
\section{Quivers}

Being a finite dimensional elementary $\CC$-algebra, $\AKSz$ can be
presented in the form $\CC Q_{n,r}/\mathcal I$, where $Q_{n,r}$ is a quiver,
$\CC Q_{n,r}$ its path algebra, and $\mathcal I$ an ideal contained in
$\mathcal J^2$ where $\mathcal J$ is the ideal generated by the edges of
$Q_{n,r}$~\cite{ARS}.
The vertices of $Q_{n,r}$ are naturally in bijection with the simple modules
$S_{[I,\c]}$, and the number $e_{[I,\c],[J,\c']}$ of edges
$S_{[I,\c]}\rightarrow S_{[J,\c']}$ is equal to

\begin{equation}
\dim \Ext^1(S_{[I,\c]},S_{[J,\c']})
=[S_{[J,\c']}  : \rad P_{[I,\c]}/\rad^2 P_{[I,\c]}].
\end{equation}
that is equal to
$ [S_{[I,\c]}  : \rad P_{[J,\c']}/\rad^2 P_{[J,\c']}]$
by auto-injectivity.

Therefore, $e_{[I,\c][J,\c']}=d_{[I,\c][J,\c']}^{(1)}$, where
\begin{equation}
d_{[I,\c][J,\c']}^{(k)}
:=[ S_{[J,\c']} : \rad^k P_{[I,\c]}/\rad^{k+1} P_{[I,\c]}]
\end{equation}
are the coefficients of the $q$-Cartan invariants
\begin{equation}
  d_{[I,\c][J,\c']}(q) =\sum_{k\ge 0} d_{[I,\c][J,\c']}^{(k)}\, q^k
\end{equation}
associated with the radical series.
Let $D_n(q)=(d_{[I,\c][J,\c']}(q))$ be the $q$-Cartan matrix.
For $n\le 4$, Section~\ref{tables} provides examples of such matrices.

As one will see in the next paragraph, the quiver of $\AKSz$ splits into
connected components spanned by evaluation classes of the colors words
of the anti-cycloribbons. In other words,
$\dim \Ext^1(S_{[I,\c]},S_{[J,\c']})\not=0$
implies that $\c$ and $\c'$ are permutations of one another.

The quivers corresponding to evaluations~$(1,3)$ and $(2,2)$ of
$\AKSzgen{4}{2}$ are given on Figures~\ref{fig:quiv13} and~\ref{fig:quiv22}.

\begin{figure}[ht]
\ifdraft




\newdimen\vcadre\vcadre=0.2cm 
\newdimen\hcadre\hcadre=0.2cm 
\def\GrTeXBox#1{\vbox{\vskip\vcadre\hbox{\hskip\hcadre%
      $\Tableau{#1\\}$%
   \hskip\hcadre}\vskip\vcadre}}

\def\arx#1[#2]{\ifcase#1 \relax \or%
  \ar @{-}[#2]  \or%
  \ar @2{-}[#2] \or%
  \ar @{--}[#2] \or%
  \ar @2{.}[#2] \or%
  \ar @{~}[#2]  \fi}

$\xymatrix@R=0.8cm@C=0.5cm{%
*{\GrTeXBox{2 & 2 & 2 & 1}}\arx1[drr] &  &
*{\GrTeXBox{2 & 2 \\ \ & 2 & 1}}\arx1[rr]&  & *{\GrTeXBox{2 \\ 2 & 2 & 1}}& &
*{\GrTeXBox{2 \\ 2 \\ 2 & 1 }}& \\
&  & *{\GrTeXBox{2 & 2 & 1 \\ \ & \ & 2}}\arx1[u]\arx1[rrd]\arx1[d]&
& *{\GrTeXBox{2 \\ 2 & 1 \\ \ & 2}}\arx1[u]\arx1[urr]\arx1[lld]\arx1[d]&  \\
&  & *{\GrTeXBox{2 & 1 \\ \ & 2 \\ \ & 2}}\arx1[d]
&  & *{\GrTeXBox{2 & 1 \\ \ & 2 & 2}}\arx1[rrd]\arx1[d]&  & \\
*{\GrTeXBox{1 \\ 2 \\ 2 \\ 2}}\arx1[rru]
&  & *{\GrTeXBox{1 \\ 2 & 2 \\ \ & 2}}
&  & *{\GrTeXBox{1 \\ 2 \\ 2 & 2}}\arx1[ll]& &
*{\GrTeXBox{1 \\ 2 & 2 & 2}}& \\
}$
\fi
\caption{\label{fig:quiv13}Quiver of $\AKSzgen{4}{2}$ corresponding to
anticycloribbons of evaluation $(1,3)$.}
\end{figure}

\begin{figure}[ht]
\ifdraft




\newdimen\vcadre\vcadre=0.2cm 
\newdimen\hcadre\hcadre=0.2cm 
\def\GrTeXBox#1{\vbox{\vskip\vcadre\hbox{\hskip\hcadre%
      $\Tableau{#1\\}$%
   \hskip\hcadre}\vskip\vcadre}}

\def\arx#1[#2]{\ifcase#1 \relax \or%
  \ar @{-}[#2]  \or%
  \ar @2{-}[#2] \or%
  \ar @{--}[#2] \or%
  \ar @2{.}[#2] \or%
  \ar @{~}[#2]  \fi}

$\xymatrix@R=0.5cm@C=0.3mm{%
*{\GrTeXBox{2 & 2 & 1 & 1} }
&  & *{\GrTeXBox{2 \\ 2 & 1 & 1} }
&  & *{\GrTeXBox{2 & 2 & 1 \\ \ &\ & 1} }
&  & *{\GrTeXBox{2 \\ 2 & 1 & \ \\ \ & 1} }& \\
&  &  & *{\GrTeXBox{2 & 1 \\ \ & 2 & 1} }\arx1[lllu]\arx1[lu]\arx1[ru]\arx1[rrru]\arx1[llld]\arx1[ld]\arx1[rd]\arx1[rrrd]& \\
*{\GrTeXBox{2 & 1 & 1 \\ \ & \ & 2}}
&  & *{\GrTeXBox{2 & 1 \\ \ & 1 \\ \ & 2} }
&  & *{\GrTeXBox{1 \\ 2 & 2 & 1} }
&  & *{\GrTeXBox{1 \\ 2 \\ 2 & 1} }& \\
&  &  & *{\GrTeXBox{1 \\ 2 & 1 \\ \ & 2} }\arx1[lllu]\arx1[lu]\arx1[ru]\arx1[rrru]\arx1[llld]\arx1[ld]\arx1[rd]\arx1[rrrd]& \\
*{\GrTeXBox{1 & 1 \\ \ & 2 & 2} }
&  & *{\GrTeXBox{1 \\ 1 \\ 2 & 2} }
&  & *{\GrTeXBox{1 & 1 \\ \ & 2 \\ \ & 2} }
&  & *{\GrTeXBox{1 \\ 1 \\ 2 \\ 2} }& \\
}$
\fi
\caption{\label{fig:quiv22}Quiver of the block $(2,2)$ of $\AKSzgen{4}{2}$.}
\end{figure}

The quiver $Q_{n,r}$ can be described for all $n$: let $[I,\c]$ and $[J,\c']$
be cycloribbons. Since the simple modules are one-dimensional, isomorphism
classes of non trivial extensions
\begin{equation}
0\rightarrow S_{[J,\c']}\rightarrow M \rightarrow S_{[I,\c]} \rightarrow 0
\end{equation}
are in one-to-one correspondence with isomorphism classes of indecomposable
two-dimensional (left) modules $M$ of socle $S_{[J,\c']}$ and such that
$M/\rad M=S_{[I,\c]}$.

Let $M$ be such a module. Since the $L_\c$ are orthogonal idempotents summing
to $1$, all but one or two of them act on $M$ by zero.

If there only is one $L_\c$ having a non-zero action on $M$, then it acts by
$1$ and the other $L_\d$ act by $0$ on $M$. So all Lagrange polynomials act as
scalars on $M$, and $\c=\c'$.
Since $M$ is an indecomposable two-dimensional module of $\AKSz$, $M$ is also
an indecomposable two-dimensional module of $H_n(0)$.
Thanks to~\cite{NCSF6}, there exist two elements $b_1$ and
$b_2$ of $M$ and an integer $i$ such that $T_i b_1=b_2$ or $(1+T_i)b_1=b_2$,
all $T_j$ for $j\in\{1,\ldots,i-2\}\cap\{i+2,\ldots,n-1\}$ act as scalars,
$T_{i-1}$ does not act as a scalar, and $T_{i+1}$ can act as a scalar or not.
Moreover, if $T_ib_1=b_2$, $T_{i-1}$ acts by $-1$ on $b_1$ and by $0$ on
$b_2$ and in the opposite way if $(1+T_i)b_1=b_2$. Finally, if $T_{i+1}$ acts
as a non-scalar, then it acts as $T_{i-1}$.

\smallskip
Since $M$ is a module, we must have $c_{i-1}=c_i=c_{i+1}$ and, it is also
equal to $c_{i+2}$ if $T_{i+1}$ acts as $T_{i-1}$. Moreover, if
$c_k\not=c_{k+1}$ then $T_k$ acts as a scalar, by $-1$ if $c_k<c_{k+1}$ and
by $0$ otherwise. One then checks that these relations are sufficient, so that
there is an edge in $Q_{n,r}$ between $[I,\c]$ and $[J,\c]$ iff there is an
edge between $I$ and $J$ in the quiver of $H_n(0)$.

Let us now consider the case where two elements, $L_\c$ and $L_{\c'}$ act as
non-zero elements on $M$. Since they commute, there exist a basis where both
are diagonal. Given that their sum is $1$ and that they are orthogonal, this
proves that there exist two elements $x$ and $y$ of $M$ such that
$x=L_\c x$ and $y=L_{\c'} y$.
Since we assumed that $M$ is indecomposable, this means that there exists an
integer $i$ such that, up to the exchange of $x$ and $y$,
$T_i x = \alpha x +\beta y$ with $\beta\not=0$.
Relations~(\ref{tx-xt}) then show that $\c'=\c\s_i$, so that $c_i\not=c_{i+1}$
and that $\alpha=0$ or $-1$, so that $T_i x = y$, or $(1+T_i)x=y$ depending on
whether $c_i>c_{i+1}$ or $c_i<c_{i+1}$.
Then, as in~\cite{NCSF6}, one proves that all $T_j$ except $T_i$ act
diagonally on $x$ and $y$, and that $T_j$ acts as a scalar if
$j\not\in\{i-1,i,i+1\}$. Finally, since $M$ is a module, the action of any
$T_j$ ($j\not=i$) on $x$ (resp. $y$) is fixed if $c_j\not=c_{j+1}$ (resp.
$c'_j\not=c'_{j+1}$). One then checks that these relations are sufficient, so
that there is an edge in $Q_{n,r}$ between $[I,\c]$ and $[J,\c']$
($\c'\not=\c$) iff there exist an integer $i$ such that $\c'=\c\s_i$, and any
integer $k\not\in \{i-1,i,i+1\}$ is a descent of both or none of $I$ and $J$.

So we can conclude:
\begin{proposition}
Let $[I,\c]$ and $[J,\c']$ be two cycloribbons.
There is an edge in $Q_{n,r}$ between $[I,\c]$ and $[J,\c]$ iff
\begin{itemize}
\item $\c=\c'$ and there is an edge in the quiver of $H_n(0)$,
\item $\c\not=\c'$ and there exists an integer $i$ such that $\c'=\c\s_i$, and
if $k$ is a descent of both or none of $I$, $J$ for $k\not\in\{i-1,i,i+1\}$.
\end{itemize}
\end{proposition}

Note that there cannot be any edge between two cycloribbons if $\c$ and $\c'$
are not permutations of one another. So $Q_{n,r}$ splits into components
generated by rearrangement classes of color vectors.
These rearrangement classes will be called \emph{evaluations} in the sequel.

Since the product of two basis elements $B_{\c,\sigma}$ and
$B_{\c',\sigma'}$ is non-zero only if $\c$ and $\c'$ are permutations of one
another, the simple modules $S_{[I,\c]}$ occurring in any indecomposable
module must be in the same evaluation class, so that the Cartan matrix of any
$\AKSz$ is a block matrix, each block being indexed by anticycloribbons with
the same evaluation.
The monochromatic blocks coincide with the Cartan matrices of $H_n(0)$. The
other one correspond to blocks (in the sense of indecomposable two-sided
ideals) of $\AKSz$.

\newpage
\section{Tables}
\label{tables}

The first tables present the $q$-Cartan matrices of the blocks $(2)$ and
$(1,1)$ of $\AKSzgen{2}{2}$ and the blocks $(3)$, $(2,1)$, $(1,2)$,
and $(1,1,1)$ of $\AKSzgen{3}{3}$. Note that the block $(1,2)$ can be
deduced from the block $(2,1)$ by changing any anticycloribbon $(I,\c)$ into
$(\overline{I},\c')$ where $\c'$ is obtained from $\c$ by exchanging its
smallest letter with its greatest, its second smallest with its second
greatest, and so on. Note also that since each permutation is the color word
of exactly one anticycloribbon, one can see a standard anticycloribbon as its
color word. Then the matrix of the block $(1,1,1)$ is the matrix whose
entry $(\sigma,\tau)$ is $q^{l(\sigma^{-1}\tau)}$.
%

Next, we present the blocks $(4)$, $(3,1)$, $(2,2)$, $(2,1,1)$, and $(1,2,1)$
of $\AKSzgen{4}{4}$.  The other blocks can be deduced from the previous ones:
$(1,1,2)$ comes from $(2,1,1)$ and $(1,1,1,1)$ is again
$q^{l(\sigma^{-1}\tau)}$.

These $q$-Cartan matrices are symmetric along both diagonals since 
the algebra is auto-injective so that the quiver is an unoriented graph (first
main diagonal) and the morphism sending $L_\c$ to $L_{\overline\c}$ and $-T_i$
to $1+T_{n-i}$ is an automorphism of $\AKSz$ (second main diagonal).
Moreover, for $q=1$, the Cartan matrix is the product of the decomposition
matrix with its transpose. 
Finally, we give three examples of decomposition matrices.

\newpage
\ifdraft

\TasseTableau

\def\tab#1{\vcenter{\vskip2pt\hbox{$\Tableau{#1}$}\vskip2pt}}
\small

\begin{figure}[ht]
\begin{center}
$\begin{array}{ccc}
 & {\tab{
1\\
1\\
}} & {\tab{
1 & 1\\
}}\\
{\tab{
1\\
1\\
}} & 1 & .\\
{\tab{
1 & 1\\
}} & . & 1\end{array}$
\qquad
$\begin{array}{ccc}
 & {\tab{
1\\
2\\
}} & {\tab{
2 & 1\\
}}\\
{\tab{
1\\
2\\
}} & 1 & q\\
{\tab{
2 & 1\\
}} & q & 1\end{array}$
\end{center}
\caption{The two blocks of the $q$-Cartan matrix of $\AKSzgen{2}{2}$.}
\end{figure}

\tiny
\begin{figure}[ht]
\begin{center}
$\begin{array}{ccccc}
 & {\tab{
1\\
1\\
1\\
}} & {\tab{
1 & \ \\
1 & 1\\
}} & {\tab{
1 & 1\\
\  & 1\\
}} & {\tab{
1 & 1 & 1\\
}}\\
{\tab{
1\\
1\\
1\\
}} & 1 & . & . & .\\
{\tab{
1 & \ \\
1 & 1\\
}} & . & 1 & q & .\\
{\tab{
1 & 1\\
\  & 1\\
}} & . & q & 1 & .\\
{\tab{
1 & 1 & 1\\
}} & . & . & . & 1\end{array}$
\qquad
\qquad
\qquad
$\begin{array}{cccccc}
 & {\tab{
1\\
1\\
2\\
}} & {\tab{
1 & 1\\
\  & 2\\
}} & {\tab{
1 & \ \\
2 & 1\\
}} & {\tab{
2 & 1\\
\  & 1\\
}} & {\tab{
2 & 1 & 1\\
}}\\
{\tab{
1\\
1\\
2\\
}} & 1 & . & q & q^2 & .\\
{\tab{
1 & 1\\
\  & 2\\
}} & . & 1 & q & . & q^2\\
{\tab{
1 & \ \\
2 & 1\\
}} & q & q & q^2 + 1 & q & q\\
{\tab{
2 & 1\\
\  & 1\\
}} & q^2 & . & q & 1 & .\\
{\tab{
2 & 1 & 1\\
}} & . & q^2 & q & . & 1\end{array}$
\end{center}
\vskip1cm
\begin{center}
$\begin{array}{cccccc}
 & {\tab{
1\\
2\\
2\\
}} & {\tab{
1 & \ \\
2 & 2\\
}} & {\tab{
2 & 1\\
\  & 2\\
}} & {\tab{
2 & \ \\
2 & 1\\
}} & {\tab{
2 & 2 & 1\\
}}\\
{\tab{
1\\
2\\
2\\
}} & 1 & . & q & q^2 & .\\
{\tab{
1 & \ \\
2 & 2\\
}} & . & 1 & q & . & q^2\\
{\tab{
2 & 1\\
\  & 2\\
}} & q & q & q^2 + 1 & q & q\\
{\tab{
2 & \ \\
2 & 1\\
}} & q^2 & . & q & 1 & .\\
{\tab{
2 & 2 & 1\\
}} & . & q^2 & q & . & 1\end{array}$
\qquad\qquad
$\begin{array}{ccccccc}
 & {\tab{
1\\
2\\
3\\
}} & {\tab{
1 & \ \\
3 & 2\\
}} & {\tab{
2 & 1\\
\  & 3\\
}} & {\tab{
3 & 1\\
\  & 2\\
}} & {\tab{
2 & \ \\
3 & 1\\
}} & {\tab{
3 & 2 & 1\\
}}\\
{\tab{
1\\
2\\
3\\
}} & 1 & q & q & q^2 & q^2 & q^3\\
{\tab{
1 & \ \\
3 & 2\\
}} & q & 1 & q^2 & q & q^3 & q^2\\
{\tab{
2 & 1\\
\  & 3\\
}} & q & q^2 & 1 & q^3 & q & q^2\\
{\tab{
3 & 1\\
\  & 2\\
}} & q^2 & q & q^3 & 1 & q^2 & q\\
{\tab{
2 & \ \\
3 & 1\\
}} & q^2 & q^3 & q & q^2 & 1 & q\\
{\tab{
3 & 2 & 1\\
}} & q^3 & q^2 & q^2 & q & q & 1\end{array}$
\end{center}
\caption{The four blocks of the $q$-Cartan matrix of $\AKSzgen{3}{3}$.}
\end{figure}

\newpage
\begin{figure}[ht]
\begin{center}
\tiny
$\begin{array}{ccccccccc}
& {\tab{
1\\
1\\
1\\
1\\
}} & {\tab{
1 & \ \\
1 & \ \\
1 & 1\\
}} & {\tab{
1 & \ \\
1 & 1\\
\  & 1\\
}} & {\tab{
1 & 1\\
\  & 1\\
\  & 1\\
}} & {\tab{
1 & 1 & 1\\
\  & \  & 1\\
}} & {\tab{
1 & 1 & \ \\
\  & 1 & 1\\
}} & {\tab{
1 & \  & \ \\
1 & 1 & 1\\
}} & {\tab{
1 & 1 & 1 & 1\\
}}\\
{\tab{
1\\
1\\
1\\
1\\
}} & 1 & . & . & . & . & . & . & .\\
{\tab{
1 & \ \\
1 & \ \\
1 & 1\\
}} & . & 1 & q & q^2 & . & . & . & .\\
{\tab{
1 & \ \\
1 & 1\\
\  & 1\\
}} & . & q & q^2 + 1 & q & . & q & . & .\\
{\tab{
1 & 1\\
\  & 1\\
\  & 1\\
}} & . & q^2 & q & 1 & . & . & . & .\\
{\tab{
1 & 1 & 1\\
\  & \  & 1\\
}} & . & . & . & . & 1 & q & q^2 & .\\
{\tab{
1 & 1 & \ \\
\  & 1 & 1\\
}} & . & . & q & . & q & q^2 + 1 & q & .\\
{\tab{
1 & \  & \ \\
1 & 1 & 1\\
}} & . & . & . & . & q^2 & q & 1 & .\\
{\tab{
1 & 1 & 1 & 1\\
}} & . & . & . & . & . & . & . & 1\end{array}$
\end{center}
\caption{The block $(4)$ of the $q$-Cartan matrix of
$\AKSzgen{4}{4}$.}
\end{figure}

\begin{figure}[ht]
\begin{center}
\tiny
$\begin{array}{ccccccccccccc}
& {\tab{
1\\
1\\
1\\
2\\
}} & {\tab{
1 & \ \\
1 & 1\\
\  & 2\\
}} & {\tab{
1 & 1\\
\  & 1\\
\  & 2\\
}} & {\tab{
1 & 1 & 1\\
\  & \  & 2\\
}} & {\tab{
1 & \ \\
1 & \ \\
2 & 1\\
}} & {\tab{
1 & 1 & \ \\
\  & 2 & 1\\
}} & {\tab{
1 & \ \\
2 & 1\\
\  & 1\\
}} & {\tab{
1 & \  & \ \\
2 & 1 & 1\\
}} & {\tab{
2 & 1\\
\  & 1\\
\  & 1\\
}} & {\tab{
2 & 1 & 1\\
\  & \  & 1\\
}} & {\tab{
2 & 1 & \ \\
\  & 1 & 1\\
}} & {\tab{
2 & 1 & 1 & 1\\
}}\\
{\tab{
1\\
1\\
1\\
2\\
}} & 1 & . & . & . & q & . & q^2 & . & q^3 & . & . & .\\
{\tab{
1 & \ \\
1 & 1\\
\  & 2\\
}} & . & 1 & q & . & q & q^2 & q^3 & q^2 & . & q^4 & q^3 & .\\
{\tab{
1 & 1\\
\  & 1\\
\  & 2\\
}} & . & q & 1 & . & q^2 & q & q^2 & q^3 & . & q^3 & q^4 & .\\
{\tab{
1 & 1 & 1\\
\  & \  & 2\\
}} & . & . & . & 1 & . & q & . & q^2 & . & . & . & q^3\\
{\tab{
1 & \ \\
1 & \ \\
2 & 1\\
}} & q & q & q^2 & . & q^2 + 1 & q^3 & q^4 + q & q & q^2 & q^3 & q^2 & .\\
{\tab{
1 & 1 & \ \\
\  & 2 & 1\\
}} & . & q^2 & q & q & q^3 & q^2 + 1 & q & q^4 + q & . & q^2 & q^3 & q^2\\
{\tab{
1 & \ \\
2 & 1\\
\  & 1\\
}} & q^2 & q^3 & q^2 & . & q^4 + q & q & q^2 + 1 & q^3 & q & q & q^2 & .\\
{\tab{
1 & \  & \ \\
2 & 1 & 1\\
}} & . & q^2 & q^3 & q^2 & q & q^4 + q & q^3 & q^2 + 1 & . & q^2 & q & q\\
{\tab{
2 & 1\\
\  & 1\\
\  & 1\\
}} & q^3 & . & . & . & q^2 & . & q & . & 1 & . & . & .\\
{\tab{
2 & 1 & 1\\
\  & \  & 1\\
}} & . & q^4 & q^3 & . & q^3 & q^2 & q & q^2 & . & 1 & q & .\\
{\tab{
2 & 1 & \ \\
\  & 1 & 1\\
}} & . & q^3 & q^4 & . & q^2 & q^3 & q^2 & q & . & q & 1 & .\\
{\tab{
2 & 1 & 1 & 1\\
}} & . & . & . & q^3 & . & q^2 & . & q & . & . & . & 1\end{array}$
\end{center}
\caption{The block $(3,1)$ of the $q$-Cartan matrix of
$\AKSzgen{4}{4}$.}
\end{figure}

\newpage
\begin{figure}[ht]
\begin{center}
\tiny
\rotateleft{
$\begin{array}{ccccccccccccccc}
 & {\tab{
1\\
1\\
2\\
2\\
}} & {\tab{
1 & \ \\
1 & \ \\
2 & 2\\
}} & {\tab{
1 & 1\\
\  & 2\\
\  & 2\\
}} & {\tab{
1 & 1 & \ \\
\  & 2 & 2\\
}} & {\tab{
1 & \ \\
2 & 1\\
\  & 2\\
}} & {\tab{
1 & \ \\
2 & \ \\
2 & 1\\
}} & {\tab{
2 & 1\\
\  & 1\\
\  & 2\\
}} & {\tab{
2 & 1 & 1\\
\  & \  & 2\\
}} & {\tab{
1 & \  & \ \\
2 & 2 & 1\\
}} & {\tab{
2 & 1 & \ \\
\  & 2 & 1\\
}} & {\tab{
2 & \ \\
2 & 1\\
\  & 1\\
}} & {\tab{
2 & 2 & 1\\
\  & \  & 1\\
}} & {\tab{
2 & \  & \ \\
2 & 1 & 1\\
}} & {\tab{
2 & 2 & 1 & 1\\
}}\\
{\tab{
1\\
1\\
2\\
2\\
}} & 1 & . & . & . & q & q^2 & q^2 & . & . & q^3 & q^4 & . & . & .\\
{\tab{
1 & \ \\
1 & \ \\
2 & 2\\
}} & . & 1 & . & . & q & . & q^2 & . & q^2 & q^3 & . & q^4 & . & .\\
{\tab{
1 & 1\\
\  & 2\\
\  & 2\\
}} & . & . & 1 & . & q & q^2 & . & q^2 & . & q^3 & . & . & q^4 & .\\
{\tab{
1 & 1 & \ \\
\  & 2 & 2\\
}} & . & . & . & 1 & q & . & . & q^2 & q^2 & q^3 & . & . & . & q^4\\
{\tab{
1 & \ \\
2 & 1\\
\  & 2\\
}} & q & q & q & q & 3\, q^2 + 1 & q^3 + q & q^3 + q & q^3 + q & q^3 + q & q^4 + 3\, q^2 & q^3 & q^3 & q^3 & q^3\\
{\tab{
1 & \ \\
2 & \ \\
2 & 1\\
}} & q^2 & . & q^2 & . & q^3 + q & q^4 + 1 & q^2 & q^2 & . & q^3 + q & q^2 & . & q^2 & .\\
{\tab{
2 & 1\\
\  & 1\\
\  & 2\\
}} & q^2 & q^2 & . & . & q^3 + q & q^2 & q^4 + 1 & . & q^2 & q^3 + q & q^2 & q^2 & . & .\\
{\tab{
2 & 1 & 1\\
\  & \  & 2\\
}} & . & . & q^2 & q^2 & q^3 + q & q^2 & . & q^4 + 1 & q^2 & q^3 + q & . & . & q^2 & q^2\\
{\tab{
1 & \  & \ \\
2 & 2 & 1\\
}} & . & q^2 & . & q^2 & q^3 + q & . & q^2 & q^2 & q^4 + 1 & q^3 + q & . & q^2 & . & q^2\\
{\tab{
2 & 1 & \ \\
\  & 2 & 1\\
}} & q^3 & q^3 & q^3 & q^3 & q^4 + 3\, q^2 & q^3 + q & q^3 + q & q^3 + q & q^3 + q & 3\, q^2 + 1 & q & q & q & q\\
{\tab{
2 & \ \\
2 & 1\\
\  & 1\\
}} & q^4 & . & . & . & q^3 & q^2 & q^2 & . & . & q & 1 & . & . & .\\
{\tab{
2 & 2 & 1\\
\  & \  & 1\\
}} & . & q^4 & . & . & q^3 & . & q^2 & . & q^2 & q & . & 1 & . & .\\
{\tab{
2 & \  & \ \\
2 & 1 & 1\\
}} & . & . & q^4 & . & q^3 & q^2 & . & q^2 & . & q & . & . & 1 & .\\
{\tab{
2 & 2 & 1 & 1\\
}} & . & . & . & q^4 & q^3 & . & . & q^2 & q^2 & q & . & . & . & 1\end{array}$
}
\end{center}
\caption{The block $(2,2)$ of the $q$-Cartan matrix of
$\AKSzgen{4}{4}$.}
\vskip3.5cm
\end{figure}

\begin{figure}[ht]
\begin{center}
\tiny
\rotateleft{
$\begin{array}{ccccccccccccccccccc}
 & {\tab{
1\\
1\\
2\\
3\\
}} & {\tab{
1 & 1\\
\  & 2\\
\  & 3\\
}} & {\tab{
1 & \ \\
1 & \ \\
3 & 2\\
}} & {\tab{
1 & 1 & \ \\
\  & 3 & 2\\
}} & {\tab{
1 & \ \\
2 & 1\\
\  & 3\\
}} & {\tab{
1 & \ \\
2 & \ \\
3 & 1\\
}} & {\tab{
1 & \ \\
3 & 1\\
\  & 2\\
}} & {\tab{
2 & 1\\
\  & 1\\
\  & 3\\
}} & {\tab{
2 & 1 & 1\\
\  & \  & 3\\
}} & {\tab{
3 & 1\\
\  & 1\\
\  & 2\\
}} & {\tab{
3 & 1 & 1\\
\  & \  & 2\\
}} & {\tab{
2 & 1 & \ \\
\  & 3 & 1\\
}} & {\tab{
1 & \  & \ \\
3 & 2 & 1\\
}} & {\tab{
3 & 1 & \ \\
\  & 2 & 1\\
}} & {\tab{
2 & \ \\
3 & 1\\
\  & 1\\
}} & {\tab{
2 & \  & \ \\
3 & 1 & 1\\
}} & {\tab{
3 & 2 & 1\\
\  & \  & 1\\
}} & {\tab{
3 & 2 & 1 & 1\\
}}\\
{\tab{
1\\
1\\
2\\
3\\
}} & 1 & . & q & . & q & q^2 & q^2 & q^2 & . & q^3 & . & q^3 & q^3 & q^4 & q^4 & . & q^5 & .\\
{\tab{
1 & 1\\
\  & 2\\
\  & 3\\
}} & . & 1 & . & q & q & q^2 & q^2 & . & q^2 & . & q^3 & q^3 & q^3 & q^4 & . & q^4 & . & q^5\\
{\tab{
1 & \ \\
1 & \ \\
3 & 2\\
}} & q & . & 1 & . & q^2 & q^3 & q & q^3 & . & q^2 & . & q^4 & q^2 & q^3 & q^5 & . & q^4 & .\\
{\tab{
1 & 1 & \ \\
\  & 3 & 2\\
}} & . & q & . & 1 & q^2 & q^3 & q & . & q^3 & . & q^2 & q^4 & q^2 & q^3 & . & q^5 & . & q^4\\
{\tab{
1 & \ \\
2 & 1\\
\  & 3\\
}} & q & q & q^2 & q^2 & q^2 + 1 & q^3 + q & 2\, q^3 & q & q & q^4 & q^4 & 2\, q^2 & q^4 + q^2 & q^5 + q^3 & q^3 & q^3 & q^4 & q^4\\
{\tab{
1 & \ \\
2 & \ \\
3 & 1\\
}} & q^2 & q^2 & q^3 & q^3 & q^3 + q & q^4 + 1 & q^4 + q^2 & q^2 & q^2 & q^3 & q^3 & q^3 + q & q^5 + q & q^4 + q^2 & q^2 & q^2 & q^3 & q^3\\
{\tab{
1 & \ \\
3 & 1\\
\  & 2\\
}} & q^2 & q^2 & q & q & 2\, q^3 & q^4 + q^2 & q^2 + 1 & q^4 & q^4 & q & q & q^5 + q^3 & q^3 + q & 2\, q^2 & q^4 & q^4 & q^3 & q^3\\
{\tab{
2 & 1\\
\  & 1\\
\  & 3\\
}} & q^2 & . & q^3 & . & q & q^2 & q^4 & 1 & . & q^5 & . & q & q^3 & q^4 & q^2 & . & q^3 & .\\
{\tab{
2 & 1 & 1\\
\  & \  & 3\\
}} & . & q^2 & . & q^3 & q & q^2 & q^4 & . & 1 & . & q^5 & q & q^3 & q^4 & . & q^2 & . & q^3\\
{\tab{
3 & 1\\
\  & 1\\
\  & 2\\
}} & q^3 & . & q^2 & . & q^4 & q^3 & q & q^5 & . & 1 & . & q^4 & q^2 & q & q^3 & . & q^2 & .\\
{\tab{
3 & 1 & 1\\
\  & \  & 2\\
}} & . & q^3 & . & q^2 & q^4 & q^3 & q & . & q^5 & . & 1 & q^4 & q^2 & q & . & q^3 & . & q^2\\
{\tab{
2 & 1 & \ \\
\  & 3 & 1\\
}} & q^3 & q^3 & q^4 & q^4 & 2\, q^2 & q^3 + q & q^5 + q^3 & q & q & q^4 & q^4 & q^2 + 1 & q^4 + q^2 & 2\, q^3 & q & q & q^2 & q^2\\
{\tab{
1 & \  & \ \\
3 & 2 & 1\\
}} & q^3 & q^3 & q^2 & q^2 & q^4 + q^2 & q^5 + q & q^3 + q & q^3 & q^3 & q^2 & q^2 & q^4 + q^2 & q^4 + 1 & q^3 + q & q^3 & q^3 & q^2 & q^2\\
{\tab{
3 & 1 & \ \\
\  & 2 & 1\\
}} & q^4 & q^4 & q^3 & q^3 & q^5 + q^3 & q^4 + q^2 & 2\, q^2 & q^4 & q^4 & q & q & 2\, q^3 & q^3 + q & q^2 + 1 & q^2 & q^2 & q & q\\
{\tab{
2 & \ \\
3 & 1\\
\  & 1\\
}} & q^4 & . & q^5 & . & q^3 & q^2 & q^4 & q^2 & . & q^3 & . & q & q^3 & q^2 & 1 & . & q & .\\
{\tab{
2 & \  & \ \\
3 & 1 & 1\\
}} & . & q^4 & . & q^5 & q^3 & q^2 & q^4 & . & q^2 & . & q^3 & q & q^3 & q^2 & . & 1 & . & q\\
{\tab{
3 & 2 & 1\\
\  & \  & 1\\
}} & q^5 & . & q^4 & . & q^4 & q^3 & q^3 & q^3 & . & q^2 & . & q^2 & q^2 & q & q & . & 1 & .\\
{\tab{
3 & 2 & 1 & 1\\
}} & . & q^5 & . & q^4 & q^4 & q^3 & q^3 & . & q^3 & . & q^2 & q^2 & q^2 & q & . & q & . & 1\end{array}$
}
\end{center}
\caption{The block $(2,1,1)$ of the $q$-Cartan matrix of
$\AKSzgen{4}{4}$.}
\end{figure}

\begin{figure}[ht]
\begin{center}
\tiny
\rotateleft{
$\begin{array}{ccccccccccccccccccc}
 & {\tab{
1\\
2\\
2\\
3\\
}} & {\tab{
1 & \ \\
2 & 2\\
\  & 3\\
}} & {\tab{
1 & \ \\
2 & \ \\
3 & 2\\
}} & {\tab{
1 & \ \\
3 & 2\\
\  & 2\\
}} & {\tab{
1 & \  & \ \\
3 & 2 & 2\\
}} & {\tab{
2 & 1\\
\  & 2\\
\  & 3\\
}} & {\tab{
2 & 1 & \ \\
\  & 3 & 2\\
}} & {\tab{
2 & \ \\
2 & 1\\
\  & 3\\
}} & {\tab{
2 & 2 & 1\\
\  & \  & 3\\
}} & {\tab{
3 & 1\\
\  & 2\\
\  & 2\\
}} & {\tab{
3 & 1 & \ \\
\  & 2 & 2\\
}} & {\tab{
2 & \ \\
3 & 1\\
\  & 2\\
}} & {\tab{
3 & 2 & 1\\
\  & \  & 2\\
}} & {\tab{
2 & \ \\
2 & \ \\
3 & 1\\
}} & {\tab{
2 & 2 & \ \\
\  & 3 & 1\\
}} & {\tab{
2 & \  & \ \\
3 & 2 & 1\\
}} & {\tab{
3 & 2 & \ \\
\  & 2 & 1\\
}} & {\tab{
3 & 2 & 2 & 1\\
}}\\
{\tab{
1\\
2\\
2\\
3\\
}} & 1 & . & q & q^2 & . & q & q^2 & q^2 & . & q^3 & . & q^3 & q^4 & q^3 & . & q^4 & q^5 & .\\
{\tab{
1 & \ \\
2 & 2\\
\  & 3\\
}} & . & 1 & q & . & q^2 & q & q^2 & . & q^2 & . & q^3 & q^3 & q^4 & . & q^3 & q^4 & . & q^5\\
{\tab{
1 & \ \\
2 & \ \\
3 & 2\\
}} & q & q & q^2 + 1 & q & q & 2\, q^2 & q^3 + q & q^3 & q^3 & q^2 & q^2 & q^4 + q^2 & 2\, q^3 & q^4 & q^4 & q^5 + q^3 & q^4 & q^4\\
{\tab{
1 & \ \\
3 & 2\\
\  & 2\\
}} & q^2 & . & q & 1 & . & q^3 & q^2 & q^4 & . & q & . & q^3 & q^2 & q^5 & . & q^4 & q^3 & .\\
{\tab{
1 & \  & \ \\
3 & 2 & 2\\
}} & . & q^2 & q & . & 1 & q^3 & q^2 & . & q^4 & . & q & q^3 & q^2 & . & q^5 & q^4 & . & q^3\\
{\tab{
2 & 1\\
\  & 2\\
\  & 3\\
}} & q & q & 2\, q^2 & q^3 & q^3 & q^2 + 1 & q^3 + q & q & q & q^4 & q^4 & q^4 + q^2 & q^5 + q^3 & q^2 & q^2 & 2\, q^3 & q^4 & q^4\\
{\tab{
2 & 1 & \ \\
\  & 3 & 2\\
}} & q^2 & q^2 & q^3 + q & q^2 & q^2 & q^3 + q & q^4 + 1 & q^2 & q^2 & q^3 & q^3 & q^5 + q & q^4 + q^2 & q^3 & q^3 & q^4 + q^2 & q^3 & q^3\\
{\tab{
2 & \ \\
2 & 1\\
\  & 3\\
}} & q^2 & . & q^3 & q^4 & . & q & q^2 & 1 & . & q^5 & . & q^3 & q^4 & q & . & q^2 & q^3 & .\\
{\tab{
2 & 2 & 1\\
\  & \  & 3\\
}} & . & q^2 & q^3 & . & q^4 & q & q^2 & . & 1 & . & q^5 & q^3 & q^4 & . & q & q^2 & . & q^3\\
{\tab{
3 & 1\\
\  & 2\\
\  & 2\\
}} & q^3 & . & q^2 & q & . & q^4 & q^3 & q^5 & . & 1 & . & q^2 & q & q^4 & . & q^3 & q^2 & .\\
{\tab{
3 & 1 & \ \\
\  & 2 & 2\\
}} & . & q^3 & q^2 & . & q & q^4 & q^3 & . & q^5 & . & 1 & q^2 & q & . & q^4 & q^3 & . & q^2\\
{\tab{
2 & \ \\
3 & 1\\
\  & 2\\
}} & q^3 & q^3 & q^4 + q^2 & q^3 & q^3 & q^4 + q^2 & q^5 + q & q^3 & q^3 & q^2 & q^2 & q^4 + 1 & q^3 + q & q^2 & q^2 & q^3 + q & q^2 & q^2\\
{\tab{
3 & 2 & 1\\
\  & \  & 2\\
}} & q^4 & q^4 & 2\, q^3 & q^2 & q^2 & q^5 + q^3 & q^4 + q^2 & q^4 & q^4 & q & q & q^3 + q & q^2 + 1 & q^3 & q^3 & 2\, q^2 & q & q\\
{\tab{
2 & \ \\
2 & \ \\
3 & 1\\
}} & q^3 & . & q^4 & q^5 & . & q^2 & q^3 & q & . & q^4 & . & q^2 & q^3 & 1 & . & q & q^2 & .\\
{\tab{
2 & 2 & \ \\
\  & 3 & 1\\
}} & . & q^3 & q^4 & . & q^5 & q^2 & q^3 & . & q & . & q^4 & q^2 & q^3 & . & 1 & q & . & q^2\\
{\tab{
2 & \  & \ \\
3 & 2 & 1\\
}} & q^4 & q^4 & q^5 + q^3 & q^4 & q^4 & 2\, q^3 & q^4 + q^2 & q^2 & q^2 & q^3 & q^3 & q^3 + q & 2\, q^2 & q & q & q^2 + 1 & q & q\\
{\tab{
3 & 2 & \ \\
\  & 2 & 1\\
}} & q^5 & . & q^4 & q^3 & . & q^4 & q^3 & q^3 & . & q^2 & . & q^2 & q & q^2 & . & q & 1 & .\\
{\tab{
3 & 2 & 2 & 1\\
}} & . & q^5 & q^4 & . & q^3 & q^4 & q^3 & . & q^3 & . & q^2 & q^2 & q & . & q^2 & q & . & 1\end{array}$
}
\end{center}
\caption{The block $(1,2,1)$ of the $q$-Cartan matrix of
$\AKSzgen{4}{4}$.}
\end{figure}

\newpage

\def\hsnew{\hfill}

\begin{figure}[ht]
$\begin{array}{ccc}
 & {\tab{
1\\
1\\
}} & {\tab{
1 & 1\\
}}\\
(2) & . & 1\\
(11) & 1 & .\end{array}$
\qquad
\qquad
\qquad
$\begin{array}{ccc}
 & {\tab{
1\\
2\\
}} & {\tab{
2 & 1\\
}}\\
(1,1) & 1 & 1\end{array}$
\caption{The two blocks $(2)$, $(1,1)$ of the decomposition matrices of
$\AKSzgen{2}{2}$.}
\end{figure}

\begin{figure}[ht]
$\begin{array}{ccccc}
 & {\tab{
1\\
1\\
1\\
}} & {\tab{
1 & \ \\
1 & 1\\
}} & {\tab{
1 & 1\\
\  & 1\\
}} & {\tab{
1 & 1 & 1\\
}}\\
(3)& . & . & . & 1\\
(21)& . & 1 & 1 & .\\
(111)& 1 & . & . & .\end{array}$
\qquad
\qquad
\qquad
$\begin{array}{cccccc}
 & {\tab{
1\\
1\\
2\\
}} & {\tab{
1 & 1\\
\  & 2\\
}} & {\tab{
1 & \ \\
2 & 1\\
}} & {\tab{
2 & 1\\
\  & 1\\
}} & {\tab{
2 & 1 & 1\\
}}\\
(2,1) & . & 1 & 1 & . & 1\\
(11,1)& 1 & . & 1 & 1 & .\end{array}$
%

\vskip1cm
$\begin{array}{cccccc}
 & {\tab{
1\\
2\\
2\\
}} & {\tab{
1 & \ \\
2 & 2\\
}} & {\tab{
2 & 1\\
\  & 2\\
}} & {\tab{
2 & \ \\
2 & 1\\
}} & {\tab{
2 & 2 & 1\\
}}\\
(1,2) & . & 1 & 1 & . & 1\\
(1,11) & 1 & . & 1 & 1 & .\end{array}$
\qquad
\qquad
\qquad
$\begin{array}{ccccccc}
 & {\tab{
1\\
2\\
3\\
}} & {\tab{
1 & \ \\
3 & 2\\
}} & {\tab{
2 & 1\\
\  & 3\\
}} & {\tab{
3 & 1\\
\  & 2\\
}} & {\tab{
2 & \ \\
3 & 1\\
}} & {\tab{
3 & 2 & 1\\
}}\\
(1,1,1) & 1 & 1 & 1 & 1 & 1 & 1\end{array}$
\caption{The blocks $(3)$, $(2,1)$, $(1,2)$, and $(1,1,1)$ of the
decomposition matrix of $\AKSzgen{3}{3}$.}
\end{figure}

\begin{figure}[ht]
\begin{center}
\rotateleft{
$\begin{array}{ccccccccc}
 & {\tab{
1\\
1\\
1\\
1\\
}} & {\tab{
1 & \ \\
1 & \ \\
1 & 1\\
}} & {\tab{
1 & \ \\
1 & 1\\
\  & 1\\
}} & {\tab{
1 & 1\\
\  & 1\\
\  & 1\\
}} & {\tab{
1 & 1 & 1\\
\  & \  & 1\\
}} & {\tab{
1 & 1 & \ \\
\  & 1 & 1\\
}} & {\tab{
1 & \  & \ \\
1 & 1 & 1\\
}} & {\tab{
1 & 1 & 1 & 1\\
}}\\
(4)   & . & . & . & . & . & . & . & 1\\
(31)  & . & . & . & . & 1 & 1 & 1 & .\\
(22)  & . & . & 1 & . & . & 1 & . & .\\
(211) & . & 1 & 1 & 1 & . & . & . & .\\
(1111)& 1 & . & . & . & . & . & . & .\end{array}$
}
\hsnew
\rotateleft{
$\begin{array}{ccccccccccccc}
 & {\tab{
1\\
1\\
1\\
2\\
}} & {\tab{
1 & \ \\
1 & 1\\
\  & 2\\
}} & {\tab{
1 & 1\\
\  & 1\\
\  & 2\\
}} & {\tab{
1 & 1 & 1\\
\  & \  & 2\\
}} & {\tab{
1 & \ \\
1 & \ \\
2 & 1\\
}} & {\tab{
1 & 1 & \ \\
\  & 2 & 1\\
}} & {\tab{
1 & \ \\
2 & 1\\
\  & 1\\
}} & {\tab{
1 & \  & \ \\
2 & 1 & 1\\
}} & {\tab{
2 & 1\\
\  & 1\\
\  & 1\\
}} & {\tab{
2 & 1 & 1\\
\  & \  & 1\\
}} & {\tab{
2 & 1 & \ \\
\  & 1 & 1\\
}} & {\tab{
2 & 1 & 1 & 1\\
}}\\
(3,1) & . & . & . & 1 & . & 1 & . & 1 & . & . & . & 1\\
(21,1) & . & 1 & 1 & . & 1 & 1 & 1 & 1 & . & 1 & 1 & .\\
(111,1) & 1 & . & . & . & 1 & . & 1 & . & 1 & . & . & .\end{array}$
}
\hsnew
\rotateleft{
$\begin{array}{ccccccccccccccc}
 & {\tab{
1\\
1\\
2\\
2\\
}} & {\tab{
1 & \ \\
1 & \ \\
2 & 2\\
}} & {\tab{
1 & 1\\
\  & 2\\
\  & 2\\
}} & {\tab{
1 & 1 & \ \\
\  & 2 & 2\\
}} & {\tab{
1 & \ \\
2 & 1\\
\  & 2\\
}} & {\tab{
1 & \ \\
2 & \ \\
2 & 1\\
}} & {\tab{
2 & 1\\
\  & 1\\
\  & 2\\
}} & {\tab{
2 & 1 & 1\\
\  & \  & 2\\
}} & {\tab{
1 & \  & \ \\
2 & 2 & 1\\
}} & {\tab{
2 & 1 & \ \\
\  & 2 & 1\\
}} & {\tab{
2 & \ \\
2 & 1\\
\  & 1\\
}} & {\tab{
2 & 2 & 1\\
\  & \  & 1\\
}} & {\tab{
2 & \  & \ \\
2 & 1 & 1\\
}} & {\tab{
2 & 2 & 1 & 1\\
}}\\
(2,2)   & . & . & . & 1 & 1 & . & . & 1 & 1 & 1 & . & . & . & 1\\
(2,11)  & . & . & 1 & . & 1 & 1 & . & 1 & . & 1 & . & . & 1 & .\\
(11,2)  & . & 1 & . & . & 1 & . & 1 & . & 1 & 1 & . & 1 & . & .\\
(11,11) & 1 & . & . & . & 1 & 1 & 1 & . & . & 1 & 1 & . & . & .\end{array}$
}
\hsnew
\rotateleft{
$\begin{array}{ccccccccccccccccccc}
 & {\tab{
1\\
1\\
2\\
3\\
}} & {\tab{
1 & 1\\
\  & 2\\
\  & 3\\
}} & {\tab{
1 & \ \\
1 & \ \\
3 & 2\\
}} & {\tab{
1 & 1 & \ \\
\  & 3 & 2\\
}} & {\tab{
1 & \ \\
2 & 1\\
\  & 3\\
}} & {\tab{
1 & \ \\
2 & \ \\
3 & 1\\
}} & {\tab{
1 & \ \\
3 & 1\\
\  & 2\\
}} & {\tab{
2 & 1\\
\  & 1\\
\  & 3\\
}} & {\tab{
2 & 1 & 1\\
\  & \  & 3\\
}} & {\tab{
3 & 1\\
\  & 1\\
\  & 2\\
}} & {\tab{
3 & 1 & 1\\
\  & \  & 2\\
}} & {\tab{
2 & 1 & \ \\
\  & 3 & 1\\
}} & {\tab{
1 & \  & \ \\
3 & 2 & 1\\
}} & {\tab{
3 & 1 & \ \\
\  & 2 & 1\\
}} & {\tab{
2 & \ \\
3 & 1\\
\  & 1\\
}} & {\tab{
2 & \  & \ \\
3 & 1 & 1\\
}} & {\tab{
3 & 2 & 1\\
\  & \  & 1\\
}} & {\tab{
3 & 2 & 1 & 1\\
}}\\
(2,1,1) & . & 1 & . & 1 & 1 & 1 & 1 & . & 1 & . & 1 & 1 & 1 & 1 & . & 1 & . & 1\\
(11,11) & 1 & . & 1 & . & 1 & 1 & 1 & 1 & . & 1 & . & 1 & 1 & 1 & 1 & . & 1 & .
\end{array}$
}
%
\hsnew
\rotateleft{
$\begin{array}{ccccccccccccccccccc}
 & {\tab{
1\\
2\\
2\\
3\\
}} & {\tab{
1 & \ \\
2 & 2\\
\  & 3\\
}} & {\tab{
1 & \ \\
2 & \ \\
3 & 2\\
}} & {\tab{
1 & \ \\
3 & 2\\
\  & 2\\
}} & {\tab{
1 & \  & \ \\
3 & 2 & 2\\
}} & {\tab{
2 & 1\\
\  & 2\\
\  & 3\\
}} & {\tab{
2 & 1 & \ \\
\  & 3 & 2\\
}} & {\tab{
2 & \ \\
2 & 1\\
\  & 3\\
}} & {\tab{
2 & 2 & 1\\
\  & \  & 3\\
}} & {\tab{
3 & 1\\
\  & 2\\
\  & 2\\
}} & {\tab{
3 & 1 & \ \\
\  & 2 & 2\\
}} & {\tab{
2 & \ \\
3 & 1\\
\  & 2\\
}} & {\tab{
3 & 2 & 1\\
\  & \  & 2\\
}} & {\tab{
2 & \ \\
2 & \ \\
3 & 1\\
}} & {\tab{
2 & 2 & \ \\
\  & 3 & 1\\
}} & {\tab{
2 & \  & \ \\
3 & 2 & 1\\
}} & {\tab{
3 & 2 & \ \\
\  & 2 & 1\\
}} & {\tab{
3 & 2 & 2 & 1\\
}}\\
(1,2,1) & . & 1 & 1 & . & 1 & 1 & 1 & . & 1 & . & 1 & 1 & 1 & . & 1 & 1 & . & 1\\
(1,11,1) & 1 & . & 1 & 1 & . & 1 & 1 & 1 & . & 1 & . & 1 & 1 & 1 & . & 1 & 1 &.
\end{array}$
}
\caption{The blocks $(4)$, $(3,1)$, $(2,2)$, $(2,1,1)$, and $(1,2,1)$ of
the decomposition matrix of $\AKSzgen{4}{4}$.}
\end{center}
\end{figure}

\newpage
\newpage
\fi

\end{document}